%Paper: alg-geom/9402008
%From: Yi Hu <yhu@math.utah.edu>
%Date: Wed, 9 Feb 94 14:02:01 EST
%Date (revised): Fri, 25 Feb 1994 15:49:34 -0500
%Date (revised): Sun, 27 Feb 1994 11:55:14 -0500
%Date (revised): Tue, 12 Dec 95 21:38:03 EST

%11/23/95

%%%%% Tex, alg-geom/9402008, revised and shorten version.

\magnification=1200

\font\smalltype=cmr10 at 10truept
 at 10truept
\pageno=0
\nopagenumbers

\def\proof{\smallskip{\sl Proof. }}
\def\endproof{\hfill\qed}
\def\qed {\nobreak$\quad$\lower 1pt\vbox{
    \hrule
    \hbox to 8pt{\vrule height 8pt\hfil\vrule height 8pt}
      \hrule}\ifmmode\relax\else\par\medbreak\fi}

\def\leftheadline{\ifnum\pageno >0
        \tenrm \folio\hfil \smalltype
         Igor Dolgachev and Yi Hu
        \hfil\fi}

\def\rightheadline{\ifnum\pageno >0
        \hfil \smalltype
 Variation of Geometric Invariant Theory Quotient
        \hfil \tenrm \folio\fi}

\headline={\ifodd\pageno\rightheadline\else\leftheadline\fi}

%$\hbox{}$
$\quad\quad\quad\quad\quad\quad\quad\quad\quad\quad\quad\quad\quad\quad$
{\smalltype To appear in Publications Math\'ematiques de l'I.H.E.S.}
\bigskip\bigskip
\centerline{\bf Variation of Geometric Invariant Theory Quotients}

%\footnote{${}^\star$}{\smalltype This Version: April 26,
%1994. Comments and suggestions will be greatly appreciated.}

\vskip .75in
\centerline{Igor V. Dolgachev  \footnote{${}^1$}{\smalltype
Research supported in part by a NSF grant.}
and Yi Hu  \footnote{${}^2$}{\smalltype
Research supported in part by  NSF grant DMS 9401695
during the revision of this paper.}
}
\bigskip
\midinsert
\narrower \noindent
{\bf Abstract.}$\quad$
 Geometric Invariant Theory gives a method for constructing quotients for group
actions on
algebraic varieties which in many cases appear as moduli spaces parametrizing
isomorphism
classes of geometric
objects (vector bundles, polarized varieties, etc.). The quotient
depends on a choice of an ample linearized line bundle.
Two choices are equivalent if they give rise to identical quotients.
A priori, there are infinitely many choices since there are infinitely many
isomorphism classes of linearized ample
line bundles. Hence several fundamental questions naturally arise. Is the set
of equivalence
classes, and hence the set of non-isomorphic quotients, finite? How does the
quotient
vary under change of the
equivalence class? In this paper we give partial answers to these questions in
the case of actions of reductive algebraic groups on nonsingular projective
algebraic varieties. We  shall show that among ample line bundles which give
projective geometric quotients there are only finitely many equivalence
classes. These classes span certain convex subsets (chambers) in a certain
convex cone in Euclidean space; and when we cross a wall separating one chamber
from another, the corresponding quotient undergoes a birational
transformation which is similar to a Mori flip.
\endinsert
\eject

\centerline{Table of contents}
\item \S 0. Introduction
\item \S 1. The numerical function $M^\bullet(x)$ and some finiteness theorems
   \item\item  \S 1.1.  The function $M^L(x)$
   \item\item  \S 1.2. Adapted one-parameter subgroups
   \item\item  \S 1.3. Stratification of the set of unstable points
\item  \S 2. Moment map and symplectic reductions
   \item\item  \S 2.1. Moment map
   \item\item  \S 2.2. Relation with geometric invariant theory
   \item\item  \S 2.3. Homological equivalence for $G$-linearized line bundles
   \item\item  \S 2.4. Stratification of the set of unstable points via moment
map
   \item\item  \S 2.5. K\"ahler quotients
\item \S 3. The $G$-ample cone
   \item\item  \S 3.1. $G$-effective line bundles
   \item\item  \S 3.2. The $G$-ample cone
   \item\item  \S 3.3. Walls and chambers
   \item\item  \S 3.4. GIT-equivalence classes
   \item\item  \S 3.5. Abundant actions
\item \S 4. Variation of quotients
   \item\item  \S 4.1. Faithful walls
   \item\item  \S 4.2. Variation of quotients
\medskip
\noindent
{\bf \S 0. Introduction}

\medskip
\noindent
{\bf 0.1} {\sl Motivation}.
Consider a projective algebraic variety $X$ acted on by a reductive algebraic
group $G$,
both defined over the field of complex numbers. In general the orbit space
 $X/G$ does not exist in the category
of algebraic varieties. One of the reasons for this is the presence of
non-closed orbits.
A way to solve this problem is suggested by
Geometric Invariant Theory. It gives a construction of a $G$-invariant open
subset $U$ of
$X$ for which the algebraic quotient $U/\!/G$ exists and $U$ is maximal with
this property.
However, there is no canonical way to choose $U$; some choice must be made and
the choice is a $G$-linearized ample line bundle.
Given an ample  $G$-linearized line bundle $L$,
one defines the set of semi-stable points as
$$X^{ss}(L) = \{x \in X: \exists \sigma  \in \Gamma (X,L^{\otimes n})^G
\quad\hbox{
such that $\sigma (x) \not= 0 \}$}$$
and the set of stable points as
$$X^s(L) =\{x\in X^{ss}(L):\quad\hbox {$G\cdot x$ is closed in $X^{ss}(L)$ and
the stabilizer $G_x$ is finite} \}.$$
As was shown by D. Mumford, a categorical quotient
$X^{ss}(L)/\!/G$ exists as a projective variety and it parametrizes the {\it
closed orbits} in $X^{ss}(L)$  and
contains the orbit space $X^s(L)/G$ as a Zariski open subset.
There are infinitely many isomorphism classes of ample  $G$-linearized line
bundles $L$.
Thereby we ask the following fundamental questions:

\item{(i)}  Is the set of non-isomorphic quotients $X^{ss}(L)/\!/G$ finite?
Descibe this set.
\item{(ii)} How does the quotient $X^{ss}(L)/\!/G$ change if we vary $L$ in the
group $Pic^G(X)$ of
isomorphism classes of $G$-linearized line bundles?

 This problem is analogous to the problem
of the variation of symplectic reductions of a symplectic manifold $M$ with
respect to an action of
a compact Lie group $K$. Recall that if $K \times M \rightarrow M$ is a
Hamiltonian action with  moment map
$\Phi: M \rightarrow {\rm Lie} (K)^{\ast}$, then for any point $p \in \Phi(M)$,
the orbit space
$\Phi^{-1}(K \cdot p)/K$ is the symplectic reduction of $M$ by $K$ with respect
to the point $p$. If
$K$ is a torus, $M=X$ with the symplectic form defined by the Chern form of
$L$, and $K$ acts on $X$
via the restriction of an algebraic action of its complexification $T$, then
the choice of a rational
point $p \in \Phi(X)$ corresponds to the choice of a $T$-linearization on $L$,
and the symplectic
reduction
$Y_p = \Phi^{-1}(p)/K$ is isomorphic to the  GIT quotient $X^{ss}(L)/\!/T$. It
turns out in this
case that if we let $p$ vary in a connected component $C$ (chamber) of the set
of regular values
of the moment map, the symplectic reductions $Y_p$ are all homeomorphic (in
fact,
diffeomorphic away from singularities) to the same manifold $Y_C$.
However if we let $p$ cross a wall separating one connected component from
another, the  reduction
$Y_p$ undergoes a very special surgery which is similar to a birational
transformation known as a flip.
This was shown in a work of V. Guillemin and S. Sternberg {\bf [GS]}. In a
purely algebraic
setting an analogous result was proved independently by M. Brion and C. Procesi
{\bf [BP]} and the second author {\bf [Hu1]}.

Our results extend the previous facts to the situation when $T$ is replaced
by any reductive group and we allow $L$ itself as well as its linearization to
vary.
We would like to point out that our  results are new even for torus actions,
considering
that we vary linearizations as well as their underlying ample line bundles in a
single
setting (the $G$-ample cone).

\medskip
\noindent
{\bf 0.2} {\sl Main Results}.
To state our partial answers to these questions, we give the following  main
definition:

\proclaim 0.2.1  Definition.
The {\it G-ample  cone}  $C^G(X)$ (for the action of $G$ on $X$) is the  convex
cone
in $NS^G(X) \otimes {\bf R}$ spanned by ample
 $G$-linearized line bundles $L$ with $X^{ss}(L)\ne\emptyset$,
where $NS^G(X)$ is the (N\'eron-Severi) group of $G$-line bundles modulo
homological  equivalence.

\smallskip
\noindent
0.2.2   One of the key ideas in our project is to introduce
certain walls in   $C^G(X)$. The philosophy is that a polarization (a
$G$-linearized
ample line bundle) lies on a wall if and only if it possesses a semi-stable
but not stable point. To this end, the Hilbert-Mumford numerical criterion of
stability
is the key clue.
 For any point $x \in X$, the Hilbert-Mumford numerical criterion of stability
allows one to introduce a function
 $Pic^G(X) \to {\bf Z}, L \to M^L(x),$ such that $x\in X^{ss}(L)$ if and only
if $M^L(x) \leq 0$,
with strict inequality for stable points.
We show that this function can be extended to a {\it lower convex} function
$M^\bullet(x)$ on $C^G(X)$.

Next we define a {\it wall} in $C^G(X)$ as the set $H(x)$ of zeroes of the
function $M^\bullet(x)$,
where the stabilizer of the point $x$ is of positive dimension.
The class $l\in C^G(X)$ of a $G$-linearized line bundle $L$
belongs to the union of walls if and only if $X^{ss}(L)\not= X^{s}(L)$. A wall
is called {\it proper} if it is not equal to the whole
$C^G(X)$. Clearly all walls are proper if and only if $X^s(L) = X^{ss}(L)$ for
some
$G$-linearized ample line bundle $L$. Only in this case are our main results
meaningful. Note that the latter condition is always satisfied when $G$ is a
torus acting faithfully.
A connected component of the complement of the union of walls is called a {\it
chamber}.
Any chamber is a convex subcone of
$C^G(X)$. Given any wall $H$, let $H^\circ$ be the set of points in $H$ which
do not lie on
any other wall unless it contains $H$ entirely. We define a {\it cell} to be
a connected component of $H^\circ$ if the latter is not empty.

We say that two $G$-linearized ample line bundles in $C^G(X)$ are {\it
$GIT$-equivalent}
if they define the same set of semi-stable points.

\proclaim 0.2.3 Theorem. Assume all walls are proper.
\item{(i)} There are only finitely many  chambers, walls and cells.
\item{(ii)} The closure of a chamber
 is a rational polyhedral cone in the interior of $C^G(X)$.
\item{(iii)} Each $GIT$-equivalence class is either a chamber or a union of
cells in the same wall.

Note that for any $G$-linearized ample line bundle representing a boundary
point
 of $C^G(X)$ the set of stable
points is empty, so they define degenerate quotients whose dimension is smaller
than the expected one.

\smallskip
\noindent
0.2.4 To answer the second main question, we describe how the  quotient changes
 when the $G$-linearized
ample line bundle
moves from one chamber to another by crossing a wall.
For special class of walls we prove that under this change the GIT quotient
undergoes
a very special birational transformation which is similar to a Mori flip.
To state the main result in this direction we need a few definitions. For any
wall $H$
 one can find a point $x\in X$
such that $H = H(x)$ and $H\ne H(y)$ for any point $y\not\in G\cdot x$ but
in the closure $\overline{G \cdot x}$ of the orbit $G\cdot x$.
Such a point is said to be pivotal for $H$. Its stabilizer is a reductive
subgroup of $G$.
A codimension 1 wall $H$ is said to be truly faithful if for any pivotal point
$x$ with $H = H(x)$
the stabilizer $G_x$ is one-dimensional (and hence isomorphic to ${\bf C}^*$).
For example, assuming that all walls are proper,
then  a codimension 1 wall is always truly faithful when $G=T$ or
when we replace $X$ by $X\times G/B$ ($G/B$ is the flag variety of $G$) and
consider
 the natural diagonal
action of $G$ on the product. The latter will assure that our theory applies to
symplectic reductions at general coadjoint orbits.
We say that two chambers
$C_1$ and $C_2$ are {\it relevant} with respect to a cell $F$ if both of them
contain $F$
in its closure and there exists a straight segment with points in
$C_1,C_2$ and in $F$.
For any chamber $C$ we denote by $X^{ss}(C)$ (resp. $X^s(C)$) the set
$X^{ss}(L)$ (resp. $X^s(L)$), where
$L$ is an ample $G$-linearized line bundle whose image is contained in $C$.
Likewise,
for any cell $F$ we denote by $X^{ss}(F)$ (resp. $X^s(F)$)
the set $X^{ss}(L)$ (resp. $X^s(L)$), where some point of $F$
can be represented by an ample $G$-linearized line bundle $L$. By Theorem 0.2.3
(iii), this definition is
independent of the choice of $L$.

\proclaim 0.2.5  Variation Theorem. Assume $X$ is nonsingular and all walls are
proper.
Let $C^+$ and $C^-$ be a pair of chambers relevant to a cell $F$ which lies on
a unique
truly faithful wall.
Then there are  two natural birational morphisms
$$f^+: X^s(C^+)/\!/G \rightarrow X^{ss}(F)/\!/G $$ and
$$f^-:  X^s(C^-)/\!/G \rightarrow X^{ss}(F)/\!/G $$
such that,
setting $\Sigma_0$ to be  $(X^{ss}(F) \setminus X^s(F))/\!/G$,  we have
\item{(i)}  $f^+$ and $f^-$ are isomorphisms over the complement of
$\Sigma_0$;
\item{(ii)} the fibres of the maps $f^\pm$ over each connected component
$\Sigma'_0 $ of  $\Sigma_0 $ are
weighted projective spaces of some dimension $d_\pm$ (the weights and
dimensions depend only on the connected component);
\item{(iii)}  $d_+ + d_-  + 1 = {\rm codim} \; \Sigma'_0$.

\medskip\noindent
{\bf 0.3}  {\sl Symplectic Reductions}.
A point in $C^G(X)$ is called a $G$-polarization on $X$. Each polarization
gives rise to a K\"ahler metric on $X$ while an integral point (polarization)
gives rise to
a Hodge metric on $X$.
In addition, for any $G$-polarization $l$, there is attached a {\it moment map}
$\Phi^l$ from $X$ to  ${\bf k}^\ast$, where
${\bf k}^\ast$ is the dual of the Lie algebra of a maximal compact subgroup $K$
of $G$.

The theory of moment maps builds a remarkable link between quotients in
Algebraic Geometry
and quotients in Symplectic Geometry.
Choose a point
$p \in \Phi(X)$, the orbit space $\Phi^{-1}(K  \cdot p)/K$ is, by definition,
the
Marsden-Weinstein symplectic reduction at $p$.
When $p$ is a regular value of $\Phi$ and the action
of $K$ on $\Phi^{-1}(K  \cdot p)$ is free, $\Phi^{-1}(K  \cdot p)/K$ inherits
the structure of a symplectic K\"ahler manifold from $X$.
It is known that for a fixed Weyl chamber $\hbar^\ast_+$,
the intersection
$\Phi(X) \cap \hbar^\ast_+$ is  a convex polytope.
The connected components of the regular values
of $\Phi$ form top chambers in $\Phi(X)\cap \hbar^\ast_+$. The set of critical
values of $\Phi$ forms walls.
If $p$ stays in a chamber, the differential structure of
$\Phi^{-1}(K  \cdot p)/K$ remains the same (the symplectic form, however, has
to change).
But when we cross a wall, the differential type
of $\Phi^{-1}(K  \cdot p)/K$ undergoes a ``flip''. This can be seen by
using the so-called {\it shifting trick}; one can
identify $\Phi^{-1}(K  \cdot p)/K$
with a K\"ahler quotient of $X \times G/B$ by the diagonal action of the group
$G$ for generic point $p$
(when $p$ is on the boundary of $\hbar^\ast_+$, one should  consider $G/P$
instead of $G/B$).
This shows that  the variation of symplectic reductions of $X$ by $K$ when $p$
is a generic
rational point in $\hbar^\ast_+$
is a special case of the variation of GIT quotients of $X \times G/B$ by $G$.
In particular our main theorems apply in this situation.

In fact, we can  introduce the so-called {\it universal moment map image}, a
region
${\cal P}$ in ${\bf k}^\ast \times C^G(X)$ which takes into account
both the K\"ahler metrics and the moment maps
$${\cal P} = \{ (v, l) \in {\bf k}^\ast \times C^G(X) : v \in \Phi^l (X) \cap
\hbar^\ast_+  \}.$$
There is a projection $\pi: {\cal P}  \rightarrow C^G(X)$ each of whose fibers
is a polytope.
This makes ${\cal P}$ a fiber-wise polytope.
Now the variational problem of symplectic reductions of $X$ by $K$ can be
treated
in this universal  moment map image (we will publish this somewhere else).

\medskip\noindent
{\bf 0.4} {\sl Hilbert, Limit Quotients, and Moduli problems}.
The quotients and the canonical morphisms among them form a projective  system.
So by taking the projective limit of this system
we obtain a variety which dominates all GIT quotients of $X$. We call this
variety the {\it limit quotient}
of $X$ by $G$.  The notion of the limit quotient is closely related to some
earlier constructions
using the Chow or
Hilbert scheme parametrizing the generic orbit closures and their limits
(see {\bf [BBS]}, {\bf [Li]},{\bf [Ka]},
{\bf [KSZ]}, and {\bf [Hu1]}).

\medskip\noindent
There are numerous examples of birational variations of moduli varieties in the
case
when the notion of stability
of geometric objects depends on a parameter. In many of these cases this
variation can be
explained as a variation of a geometric invariant quotient.
This basic observation was probably made first by M. Thaddeus (cf. {\bf
[Re],[Th1]}).
After the first preprint version of our paper had appeared, we have received a
preprint of
Thaddeus (now the paper {\bf [Th2]}) where some of our results were reproven by
different and original methods. For example, he finds the structure of flips by
using Luna's Slice Theorem.
This allows him to obtain some further information about the algebraic
structure of the flip maps in Theorem 0.2.5. For instance, he proves that the
weighted projective space fibrations of the exceptional loci are locally
trivial.
A similar approach to the proof of Theoren 0.2.5 was independently
proposed by Charles Walter.
In this paper we are not discussing the applications of our theory to the limit
quotients
and  the moduli problems,  and plan to
return to this and produce more examples in subsequent publications.

\medskip\noindent
{\bf Acknowledgements}.  We would like to thank the following people who helped
us
in various ways during the preparation of this work: Steve Bradlow, Michel
Brion, J\'anos Koll\'ar,
Miles Reid, Michael Thaddeus, among others.
We are especially grateful to Charles Walter who suggested
counter examples to the earlier
versions of our Variation Theorem.
Thanks are also due to the referee for making many useful critical comments.
Finally we happily acknowledge the
hospitality of the  Max-Planck-Institut where part of this work was written up.

\bigskip

\noindent
{\bf \S1. The numerical function $M^\bullet (\bullet)$ and some finiteness
theorems}
% }
\smallskip
Throughout the paper we
freely use the terminology and the basic facts
from Geometric Invariant Theory which
can be found in {\bf [MFK]}. We shall
be working over the field ${\bf C}$ of complex numbers, although
most of what follows is valid over an
arbitrary algebraically closed field.

\smallskip
\noindent  {\bf 1.1} {\sl The function $M^L(x)$ \/}.
\smallskip\noindent
1.1.1 Let
$\sigma: G\times X\rightarrow X$ be an algebraic action of a connected
reductive linear
algebraic group $G$ on an irreducible projective algebraic
variety $X$. Let $L\in Pic^G(X)$ be a $G$-linearized line bundle over $X$. For
every
$x\in X$ and any 1-parameter subgroup $\lambda:{\bf C^*}\to G$,  the
1-parameter subgroup
$\lambda ({\bf C^*})$ acts on the fiber $L_{x_0}$ via the charactor
$t^{\mu^L(x,\lambda)}$
where $x_0 = \lim_{t\to 0}{\lambda}(t)\cdot x$. The number
$\mu^L(x,\lambda)$  satisfies the following properties:
\item{i)} $\mu^L(g\cdot x,g\cdot \lambda\cdot g^{-1})=\mu^L(x,\lambda)$ for any
$g\in G$;
\item{ii)} for fixed $x$ and $\lambda$, the map $L\mapsto \mu^L(x,\lambda)$ is
a
homomorphism $Pic^G(X) \to {\bf Z}$;
\item{iii)} $\mu^L(x,g\cdot \lambda\cdot g^{-1}) =\mu^L (x,\lambda)$ for any
$g\in P(\lambda),$
where $P(\lambda)$ is a certain parabolic subgroup
of $G$ associated with $\lambda$ (consult 1.2.1 in the sequel);
\item{iv)} $\mu^L(\lim_{t\to 0}{\lambda}(t)\cdot x,\lambda) =
\mu^L(x,\lambda)$.

\smallskip
The numbers $\mu^L(x,\lambda)$ are used to give a numerical criterion for
stability of
points in $X$ with respect to an ample $G$-linearized line bundle $L$:
$$ x\in X^{ss}(L) \Leftrightarrow \mu^L(x,\lambda) \le 0 \quad \hbox {for all
one-parameter
subgroups
$\lambda$},$$
$$ x\in X^{s}(L) \Leftrightarrow \mu^L(x,\lambda) < 0 \quad\hbox {for all
one-parameter
subgroups
$\lambda$}.$$
As usual, $X^{us}(L):=X \setminus X^{ss}(L)$ denotes the set of unstable
points.
For future use we denote the set  $X^{ss}(L) \setminus X^{s}(L)$ by
$X^{sss}(L)$. Its
points are called {\it strictly semistable}.

\smallskip\noindent
1.1.2
It follows from the numerical criterion of stability that
$$\displaylines{X^{ss}(L) = \bigcap_{T\  {\rm maximal\  torus}}X^{ss}(L_T),\cr
X^{s}(L) = \bigcap_{T\  {\rm maximal\  torus}}X^{s}(L_T),\cr}$$
where $L_T$ denotes the image of $L$ under the restriction map $Pic^G(X)\to
Pic^T(X)$.

\smallskip\noindent
1.1.3 Let $T$ be a maximal torus of $G$ and $W = N_G(T)/T$ be its Weyl group.
We denote by ${\cal X}_*(G)$ the set of one-parameter subgroups of $G$. We have
$$
{\cal X}_*(G) = \bigcup_{g\in G}{\cal X}_*(gTg^{-1}).
$$
Let us identify ${\cal X}_*(T)\otimes {\bf R}$ with  ${\bf R}^n$  and consider
a $W$-invariant
Euclidean norm $\|\ \|$ in ${\bf R}^n$. Then we can define for any
$\lambda \in {\cal X}_*(G)$,
$$
\|\lambda \|: = \|Int(g)\circ \lambda \|,
$$
where $Int(g)$ is the inner automorphism of $G$ such that $Int(g)\circ \lambda
\in {\cal X}_*(T)$.

Now let
$$\displaylines{
\bar{\mu}^L(x,\lambda):= {\mu^L(x,\lambda) \over \|\lambda\|}, \cr
M^L(x) = sup_{\lambda \in {\cal X_*}(G)} \bar{\mu}^L(x,\lambda).
\cr}$$
We shall show in Proposition 1.1.6 that $M^L(x)$ is always finite (cf. also
1.1.4 below).
This function plays a key role in the rest of the paper.

\proclaim 1.1.4 Lemma. Assume $L$ is ample.
 Let $T$ be a maximal torus of $G$ and $r_T:Pic^G(X) \to Pic^T(X)$ be the
restriction map.
Then for any $x\in X$,
the set $\{M^{r_T(L)}(g\cdot x), g\in G\}$ is finite and $M^L(x) = max_{g\in
G}M^{r_T(L)}(g\cdot x)$.

\proof  See {\bf [Ne1]}, Lemma 3.4.
\endproof

\noindent
1.1.5 Assume $L$ is ample. One can give the following interpretation of the
function $M^L(x)$.
Replacing $L$ with some positive power $L^{\otimes n}$,
 we may assume that $L$ is very ample.
 Choose a $G$-equivariant embedding of $X$ in a projective space ${\bf P}(V)$
such that $G$ acts on $X$ via its linear representation in $V$.
Let us assume that $G$ is the $n$-dimensional
torus $({\bf C}^*)^n$. Then its group of characters ${\cal X}(G)$ is isomorphic
to
${\bf Z}^n$. The isomorphism is defined by assigning to any
$(m_1,\ldots,m_n) \in {\bf Z}^n$ the homomorphism $\chi:G\to {\bf C}^*$ defined
by the formula
$$
\chi((t_1,\ldots,t_n)) = t_1^{m_1}\cdot \ldots \cdot t_n^{m_n}.
$$
Any one-parameter subgroup ${\lambda}:{\bf C}^* \to G$ is given by the formula
$$
{\lambda}(t) = (t^{r_1},\ldots,t^{r_n})
$$
for some $(r_1,\ldots,r_n)\in {\bf Z}^n$. In this way we can identify the set
${\cal X}_*(G)$
 with the
 group ${\bf Z}^n$. Let $\lambda \in {\cal X}_*(G)$ and  $\chi \in {\cal
X}(G)$;
the composition ${\chi}\circ {\lambda}$ is a homomorphism ${\bf C}^* \to {\bf
C}^*$, hence is
defined by an integer $m$. We denote this integer by $\langle \lambda,\chi
\rangle$.
It is clear that the pairing
$$
{\cal X}_*(G)\times {\cal X}(G) \to {\bf Z},
$$
$$ (\lambda,\chi) \mapsto {\langle \lambda,\chi \rangle}$$
is isomorphic to the natural dot-product pairing
$$
{\bf Z}^n \times {\bf Z}^n \to {\bf Z}.
$$
Now let
$$
V = \bigoplus_{\chi \in {\cal X}(G)}V_{\chi},
$$
where $V_{\chi} = \{v\in V: g\cdot v = {\chi}(g)\cdot v\}$. For any $v\in V$ we
can write
$v = \sum_{\chi}v_{\chi}$, where $v_{\chi} \in V_{\chi}$.
The group $G$ acts on a vector $v$ by the formula
$$
	g\cdot v = \sum_{\chi}{\chi}(g)\cdot v.
$$
Let $x\in {\bf P}(V)$ be represented by a vector $v$ in $V$. We set
$$\openup2pt \displaylines{
st(x) = \{\chi:v_{\chi} \not= 0\} \quad\hbox{(the state set of $x$)},\cr
\overline {st(x)} = \quad\hbox {convex hull of st($x$) in ${\cal
X}^*$}(G)\otimes {\bf R} \cong {\bf R}^n.\cr}
$$
Then
$$
\mu^L(x,\lambda) = \min_{\chi\in st(x)} \langle \lambda,\chi \rangle.
$$
In particular,
$$\openup2pt \displaylines{
x\in X^{ss}(L) \Leftrightarrow 0 \in \overline{st(x)} \cr
x\in X^s(L) \Leftrightarrow 0 \in int\{\overline{st(x)}\}.\cr}
$$

In fact,
${\mu^L(x,\lambda) \over \|\lambda\|}$ is
equal to the signed distance from the origin to the boundary of
the projection of $\overline {st(x)}$ to the positive ray spanned by the vector
$\lambda$.
Then $|M^L(x)|$ is equal to the distance from the origin to the boundary of
$\overline {st(x)}$.
Now if $G$ is any reductive group we can fix a maximal torus $T$ in $G$ and
apply Lemma 1.1.4.
This will give us the interpretation
of the function $M^L(x)$ in the general case.

\proclaim 1.1.6 Proposition.  For any $L\in Pic^G(X), M^L(x)$ is finite.

\proof  It follows from the previous discussion that $M^L(x)$ is finite if $L$
is ample.
It is known that for
any  $L\in Pic(X)$ and an ample $L_1 \in Pic(X)$ the bundle $L\otimes
L_1^{\otimes N}$
is ample for sufficiently large $N$. This shows that any $L\in Pic^G(X)$ can be
written as a difference
$L_1\otimes L_2^{-1}$ for ample $L_1,L_2\in Pic^G(X)$. We have
$$M^{L_2^{-1}}(x) = \sup_\lambda {\mu^{L_2^{-1}}(x,\lambda) \over \|\lambda\|}
=
\sup_\lambda {-\mu^{L_2}(x,\lambda) \over \|\lambda\|} = -\inf_\lambda
{\mu^{L_2}(x,\lambda) \over \|\lambda\|}.$$

If $G$ is a torus, then it follows from the last paragraph of 1.1.5
that $\inf_\lambda {\mu^{L_2}(x,\lambda) \over \|\lambda\|}$ is finite.
One  can also argue  as follows.
The function ${\lambda \over \|\lambda\|} \to {\mu^{L_2}(x,\lambda)
\over \|\lambda\|}$ extends to a continuous function
on the sphere of radius 1. Hence it is bounded from below and from above. If
$G$ is any reductive group,
 and $T$ its maximal torus,
we use that ${\mu^{L_2}(x,\lambda) \over \|\lambda\|} = {\mu^{r_T (L_2)}(g\cdot
x,\lambda_1) \over \|\lambda_1\|}$
for some $g\in G$ and
some $\lambda_1\in {\cal X}_*(T)$. Since the set of possible state sets
$st_T(g\cdot x)$ is finite, we obtain that
$$\inf_\lambda {\mu^{L_2}(x,\lambda) \over \|\lambda\|} =
\inf_g\inf_{\lambda_1}{\mu^{r_T(L_2)}(g\cdot x,\lambda_1) \over
\|\lambda_1\|}$$
is finite. Now
$$M^L(x) = \sup_\lambda ({\mu^{L_1}(x,\lambda) \over
\|\lambda\|}+{\mu^{L_2^{-1}}(x,\lambda) \over \|\lambda\|}) \leq
\sup_\lambda {\mu^{L_1}(x,\lambda) \over \|\lambda\|}+\sup_\lambda
{\mu^{L_2^{-1}}(x,\lambda) \over \|\lambda\|}$$
which implies that $M^L(x)$ is finite because $M^L(x)$ is obviously bounded
from below.
\endproof

\noindent
1.1.7 One can restate the numerical criterion using the function $M^L(x)$:
For any ample $L\in Pic^G(X)$ on a complete X
$$\displaylines{
X^{ss}(L) = \{x\in X: M^L(x) \le 0\},\cr
X^{s}(L) = \{x\in X: M^L(x) < 0\}.\cr}
$$

\medskip\noindent
{\bf 1.2} {\sl Adapted one-parameter subgroups}. The main
references here are {\bf [Ne2], [Ke], [Ki]}.

\smallskip\noindent
1.2.1 For every one-parameter subgroup $\lambda$ one defines a subgroup
$P(\lambda) \subseteq G$ by
$$
P(\lambda):= \{g\in G: \lim_{t\to 0}\lambda (t)\cdot g \cdot \lambda (t)^{-1}
\quad\hbox {exists in
$G$}\}.
$$

This is the parabolic subgroup of $G$; $\lambda$ is contained in its radical.
The set
$$ L(\lambda): = \{ \bar{g} =
\lim_{t\to 0}\lambda (t)\cdot g\cdot \lambda(t)^{-1} : g\in P(\lambda) \}
$$
is the subgroup of
$P(\lambda)$ which centralizes $\lambda$, the set of $g$'s such that the limit
equals $1$ forms
the  unipotent
radical $U(\lambda)$ of $P(\lambda)$ (see {\bf [MFK]}, Prop. 2.6). We have
$$P(\lambda) = U(\lambda)  L(\lambda) \quad  \hbox{(semi-product)}.$$
Let
$$ Lie(G)  = {\bf t} \oplus \bigoplus_{\alpha\in \Phi} Lie(G)_{\alpha}$$
be the root decomposition for the Lie algebra of $G$, where ${\bf t}=$
Lie($T$). Then
$$  Lie(L(\lambda)) = {\bf t} \oplus \bigoplus_{\langle \lambda,\alpha \rangle
= 0} Lie(G)_{\alpha},
\ \ Lie(U(\lambda)) = \bigoplus_{\langle \lambda,\alpha \rangle > 0}
Lie(G)_{\alpha}.$$

\proclaim 1.2.2 Definition. Let $L\in Pic^G(X), x\in X$. A one-parameter
 subgroup $\lambda$ is called {\it adapted} to $x$ with respect to $L$ if
$$M^L(x) = {\mu^L(x,\lambda)\over \|\lambda\|}.$$
The set of primitive (i.e. not divisible by any positive integer) adapted
one-parameter subgroups will be denoted by $\Lambda^L(x)$.

\proclaim 1.2.3 Corollary. Assume $L$ is ample. The set $\Lambda^L(x)\cap {\cal
X}_*(T)$ is non-empty
and finite.

\proof
 This follows from observing that the set of all possible states $st(x)$ with
respect to $r_T(L)$ is finite.
\endproof

\proclaim 1.2.4 Theorem. Assume $L$ is ample and $x\in X^{us}(L)$. Then
\item {(i)} There exists a parabolic subgroup $P(L)_x \subseteq G$ such that
$P(L)_x = P(\lambda)$
for all $\lambda \in \Lambda(L)_x$.
\item {(ii)} All elements of $\Lambda^L(x)$ are conjugate to each other
by elements from $P(\lambda)$.
\item {(iii)} If $T$ is a maximal torus contained in $P(L)_x$, then
$\Lambda^L(x)\cap {\cal X}_*(T)$
forms one orbit with respect to the action of the Weyl group.

\proof  This is a theorem of G. Kempf ${\bf [Ke]}$ .
\endproof

\proclaim 1.2.5 Corollary. Assume $x\in X^{us}(L)$. Then for all $g\in G$
\item{(i)} $\Lambda (gx) = Int(g) \Lambda (x)$;
\item{(ii)} $P(L)_{g\cdot x} = Int(g) P(L)_x $;
\item{(iii)} $G_x\subset P(L)_x$.

\proclaim 1.2.6 Theorem. Assume $L\in Pic^G(X)$ is ample and  $x \in
X^{us}(L)$.  Let $\lambda \in {\Lambda}(L)_x$
and $y=\lim_{t\to 0}{\lambda}(t)\cdot x$. Then
\item{(i)} $\lambda \in \Lambda (L)_y$;
\item{(ii)} $M^L(x) = M^L(y)$;

\proof This is Theorem 9.3 from {\bf [Ne2]}.
\endproof

\medskip\noindent
{\bf 1.3}. {\sl Stratification of the set of unstable points}.
Following Hesselink {\bf [He]\/},  we introduce the following stratification
of the set $X^{us}(L)$.

\smallskip\noindent
1.3.1  For each $d > 0$ and conjugacy class
$<\tau>$ of a one-parameter subgroup $\tau$ of $G$, we set
$$S^L_{d,<\tau>} = \{x\in X:M^L(x) = d,\exists g\in G \quad\hbox{such that}\
Int(g)\circ \tau \in \Lambda(L)_x\}.$$
Let ${\cal T}$ be the set of conjugacy classes of one-parameter subgroups.
Hesselink shows that for any ample $L$
$$ X=X^{ss}(L)\cup \bigcup_{d>0,<\tau>\in {\cal T}}S^L_{d,<\tau>}$$
is a finite stratification of $X$ into Zariski locally closed $G$-invariant
subvarieties of $X$.
Observe that property (ii) of Theorem 1.2.4 ensures that the subsets
$S_{d,<\tau>}, d > 0,$ are disjoint.

\smallskip\noindent
1.3.2 Let $x\in S^L_{d,<\tau>}$. For each $\lambda \in \Lambda(x)$ we can
define the point
$y=\lim_{t\to 0}{\lambda}(t)\cdot x$. By Theorem 1.2.6, $y\in S^L_{d,<\tau>}$
and $\lambda$ fixes $y$.
Set
$$Z^L_{d,<\tau>}:= \{x\in S^L_{d,<\tau>}:{\lambda}({\bf C}^*) \subset G_x
\quad\hbox{ for some}\  \lambda \in <\tau>\}.$$
Let us further subdivide each stratum $S^L_{d,<\tau>}$ by putting for each
${\lambda}\in <\tau>$
$$ S^L_{d,\lambda}:= \{x\in S^L_{d,<\tau>}:\lambda \in \Lambda^L (x)\}.$$
Hesselink calls these subsets {\it blades}. If
$$Z^L_{d,\lambda}:= \{x\in S^L_{d,\lambda}:\lambda({\bf C}^*) \subset G_x\},$$
then we have the map
$$p_{d,\lambda}: S^L_{d,\lambda}\to Z^L_{d,\lambda}, \; x\mapsto y=\lim_{t\to
0}\lambda(t)\cdot x.$$

\proclaim 1.3.3 Proposition.
\item {(i)}$S^L_{d,\lambda} = \{x\in X:\lim_{t\to 0}\lambda(t)\cdot x \in
Z^L_{d,\lambda}\} =p_{d,\lambda}^{-1} (Z^L_{d,\lambda})$;
\item {(ii)} for each connected component $Z^L_{d,\lambda ,i}$ of
$Z^L_{d,\lambda}$ the restriction of the map
$p_{d,\lambda}:S^L_{d,\lambda}\to Z^L_{d,\lambda}$ over $Z^L_{d,\lambda,i}$ is
a vector bundle with the zero section
equal to $Z^L_{d,\lambda,i}$, assuming in addition that $X$ is smooth;
\item{(iii)}$S^L_{d,\lambda}$ is $P(\lambda)$ -invariant, $Z^L_{d,\lambda ,i}$
is  $L(\lambda)$-invariant;
 if $\bar g$ denotes the projection of $g\in P(\lambda)$ to $L(\lambda)$, then
for each $x\in S^L_{d,\lambda}, p_{d,\lambda}(g\cdot x) = \bar g\cdot
p_{d,\lambda}(x)$;
\item {(iv)} there is a surjective finite morphism $
G\times_{P(\lambda)}S^L_{d,\lambda} \to S^L_{d,<\tau>}$ .
It is bijective if $d > 0$ and is an isomorphism if
$S^L_{d,<\tau>}$ is normal.

\proof See {\bf [Ki]}, \S13.2, 13.5.
\endproof
\smallskip\noindent
1.3.4 Let $X^\lambda = \{x\in X:\lambda({\bf C}^*)\subset G_x\}.$ By
definition, $Z^L_{d,\lambda}\subset X^\lambda$.
As usual we may assume that $G$ acts on $X$ via its linear representation in
the space $V = \Gamma (X,L)^*$.
Let $V = \oplus_iV_i$, where $V_i = \{v\in V: \lambda (t)\cdot v = t^iv\}$.
Then $X^\lambda = \cup_iX^\lambda_i$,
where $X^\lambda_i = {\bf P}(V_i)\cap X$. For any $x\in Z^L_{d,\lambda}, d
=M^L(x) = \mu^L(x,\lambda)$. If $v$ represents
$x$ in $V$, then, by definition of $\mu^L(x,\lambda)$, we have $v\in V_d$.
Therefore, $Z^L_{d,\lambda}\subset X^\lambda_d$.
Since $L(\lambda)$ centralizes $\lambda$ each $\lambda$-eigensubspace $V_i$
is stable with respect to $L(\lambda)$. By Proposition 1.3.3 (iii), the group
$L(\lambda)$ leaves $Z^L_{d,\lambda}$ invariant. The central subgroup
$\lambda({\bf C}^*)$
of $L(\lambda)$ acts identically on ${\bf P}(V_d)$ , hence one defines the
action of
$L(\lambda)' = L(\lambda)/\lambda({\bf C}^*)$ on ${\bf P}(V_d)$ leaving
$Z^L_{d,\lambda}$ invariant.
Let ${\cal O}_{X^\lambda_d}(1) $ be the very ample $L(\lambda)'$-linearized
line bundle obtained from the embedding
of $X^\lambda_d$ into ${\bf P}(V_d)$.

\proclaim 1.3.5 Proposition.  $Z^L_{d,\lambda} = (X^\lambda_d)^{ss}({\cal
O}_{X^\lambda_d}(1))$.

\proof See {\bf [Ne2]},Theorem 9.4.
\endproof

\proclaim 1.3.6 Lemma. Let $\Pi(X^\lambda)$ be the set of connected components
of $X^\lambda,
\lambda \in {\cal X}_*(G)$. Then the
 group $G$ acts on the union of the sets $\Pi(X^\lambda)$ ($\lambda \in {\cal
X}_*(G)$)
 via its action on the set ${\cal X}_*(G)$ with finitely many orbits.

\proof
 Choose some $G$-equivariant projective embedding of $X \hookrightarrow {\bf
P}(V)$.
Let $T$ be a maximal torus of $G$ and $V = \oplus_\chi V_\chi$ be the
decomposition
into the direct sum of eigenspaces with respect to the action of $T$.
Since any $\lambda$ is conjugate to some one-parameter subgroup of $T$, it is
enough to show that
the set $(\cup_{\lambda\in {\cal X}_*(T)}\Pi(X^\lambda))$ is a finite union.
For each integer $i$
and $\lambda$, let $V_{i, \lambda} = \oplus_{\langle \chi,\lambda \rangle = i}
V_\chi$. Then
$X^\lambda = \cup_i{\bf P}(V_{i, \lambda}) \cap X$. Since the number of
non-trivial direct summands
$V_\chi$ of $V$ is finite, there are only finitely many different subvarieties
of
the form $X^\lambda$. Each of them consists of finitely many
connected components.
\endproof

\smallskip\noindent
1.3.7 For each connected component $X^{\lambda}_i$ of $X^{\lambda}$, one can
define its contracting set
$X^+_i = \{x\in X:\lim_{t\to 0}\lambda(t)\cdot x \in X^{\lambda}_i\}$. By a
theorem of Bialynicki-Birula {\bf [B-B2]},
this set has the structure of a vector bundle with respect to the natural map
$x\mapsto \lim_{t\to 0}\lambda(t)\cdot x$
when $X$ is smooth. The decomposition $X = \cup_i X^+_i$ is the so-called
Bialynicki-Birula decomposition induced by $\lambda$.
Although there are infinitely many one-parameter subgroups of $G$,
the number of their corresponding Bialynicki-Birula decompositions
is finite up to the action of the group $G$. This follows from the next
proposition which is proven in {\bf [Hu2]}:

\proclaim 1.3.8 Proposition. Let $T$ be a fixed maximal torus of $G$. Then
there are only finitely many
Bialynicki-Birula decompositions induced by $\lambda \in {\cal X}_*(T)$.

\smallskip
In fact, there is a decomposition of ${\cal X}_*(T)$ into a finite union of
rational cones such that
two  one-parameter subgroups give rise to the same Bialynicki-Birula
decomposition if and only if they lie in the same cone.
This is Theorem 3.5 of {\bf [Hu2]}. A simple example of this theorem is the
case of the flag variety $G/B$
acted on by a maximal torus $T$. In this case,
the cone decompostion of  ${\cal X}_*(T)$ is just the fan formed by the Weyl
chambers. Each Weyl chamber gives
a Bialynicki-Birula decomposition (there are $|W|$ of them, where $W$ is the
Weyl group).
 Other Bialynicki-Birula decompositions correspond to the faces of the
Weyl chambers.

\smallskip\noindent
One can also prove this proposition without appealing to the cone decomposition
of ${\cal X}_*(T)$.
 Let $W = \cup_i  W_i$ be the union of the connected components
of the fixed point set of $T$. Two points $x$ and $y$ of $X$ are called
equivalent if whenever $\overline{T \cdot x} \cap W_i
\ne \emptyset$, then $\overline{T \cdot y} \cap W_i
\ne \emptyset$, and vice versa.  This gives a decomposition of $X$ into a
finite union of equivalence classes
$X = \cup_E X^E$,
 where $E$ ranges over all equivalence classes ($X^E$ are called torus strata
in {\bf [Hu2]}).
 It is proved in Lemma 6.1 of {\bf [Hu2]}
that every Bialynicki-Birula stratum $X^+_i$ is a (finite) union of $X^E$. This
implies that
there are only finitely many  Bialynicki-Birula decompositions.

\proclaim 1.3.9 Theorem.
\item {(i)} The set of locally closed subvarieties $S$ of $X$ which can be
realized as the stratum
$S^L_{d,<\tau>}$ for some ample $L\in Pic^G(X),d>0$ and $\tau \in {\cal
X}_*(G)$ is finite.
\item {(ii)} The set of possible open subsets of $X$ which can be realized as a
subset of semi-stable
 points with respect to some ample $G$-linearized line bundle is finite.

\proof The second assertion obviously follows from the first one.
We prove the first assertion
by using induction on the rank of $G$.
The assertion is obvious if $G = \{1\}$. Let us apply the induction to
the case when $X$ is equal to a connected component $X_i^\lambda$ of
$X^\lambda$ ($\lambda \in {\cal X}_*(T)$)
and $G = L(\lambda)/T_i$ where $T_i$ is the isotropy of $X_i^\lambda$ in $T$.
It is well known that only finitely many
subgroups of $T$ can be realized as isotropy subgroups.
Thus we obtain that  the set
of open subsets of $X_i^\lambda$ which can be realized as the sets of
semi-stable points of some ample
 $L(\lambda)/T_i$-linearized line bundle
is finite. By Proposition 1.3.5 and Lemma 1.3.6,
 this implies that the set of locally closed subsets of $X$ which can be
realized as the subsets
$Z^L_{d,<\tau>}$ is finite. Now, by Proposition 1.3.3 (i) and Proposition
1.3.8,
out of the  finitely many $Z^L_{d,<\tau>}$, one can only construct finitely
many $S^L_{d,<\tau>}$, and we are done.
\endproof

\bigskip\noindent
{\bf \S 2.  Moment map and symplectic reductions.}

\smallskip
Here we explain the relationship between the geometric invariant theory
quotients
(briefly GIT quotients) and the symplectic reductions.
The main references are {\bf [Au], [Ki], [Ne2]}. The discussion is necessary
for extending
 geometric invariant theory to the K\"ahler case, which is useful throughout \S
3.

\smallskip\noindent
{\bf 2.1} {\sl Moment map}.
\smallskip\noindent
2.1.1  Let $M$
 be a compact symplectic manifold, i.e. a compact smooth manifold equipped with
a non-degenerate
smooth closed 2-form ${\omega}$ on $M$ (a symplectic form). The non-degeneracy
condition means that for any
point $x\in M$ the skew-symmetric bilinear form $\omega_x$ on the tangent space
$T(M)_x$ is non-degenerate, i.e. defines an isomorphism
$T(M)_x \to T(M)_x^*$. This implies in particular that $M$
 is even-dimensional. The most relavant example for us is when $M$ is a
nonsingular closed algebraic
subvariety of  projective space equipped with the Fubini-Study symplectic form,
the imaginary part of
the standard K\"ahler form of projective space.

Let $K$ be a compact
Lie group which acts {\it symplectically} on a symplectic manifold
$(M,\omega)$. This means that $K$ acts smoothly on $M$ and
for any $g\in K, g^*(\omega) =\omega$.
We denote by $Lie(K)$ the Lie algebra of $K$. We shall consider any $\xi\in
Lie(K)$ as a linear function on the dual space $Lie(K)^*$.
Each point $x\in M$ defines a map $K\to M, g\mapsto g\cdot x$. For any $\xi\in
Lie(K)$, the differential of this map at the identity element
of $K$ sends $\xi$ to a tangent vector $\xi^\sharp_x \in T(K\cdot x)_x \subset
T(M)_x$. Thus each $\xi$ in $Lie(K)$ defines a vector field
$\xi^\sharp\in T(M), x\mapsto \xi^\sharp_x$. The non-degeneracy condition
allows one to define an isomorphism $\imath_{\omega}$ from the space $T(M)$ of
smooth vector fields on $M$
to the space $T^*(M)$ of smooth 1-forms on $M$.

\proclaim 2.1.2 Definition. A {\it moment map} for the action of $K$ on $M$ is
a smooth map
$\Phi:M\to Lie(K)^*$ satisfying the following two properties:
\item {(i)} $\Phi$ is equivariant with respect to the action of $K$ on $M$ and
the
co-adjoint action of $K$ on $Lie(K)^*$,
\item {(ii)} for any $\xi\in Lie(K)$, $\imath_{\omega}(\xi^\sharp) = d(\xi\circ
\Phi)$.

\smallskip
A moment map always exists if $H^1(M,{\bf R}) = 0$. It will always exist in our
situation. It
is defined uniquely up to the addition of
an element from $Lie(K)^*$ which is invariant with respect to the coadjoint
action. In particular, it is unique if $K$ is semi-simple.

 The next result describes the derivative of the moment map.

\proclaim 2.1.3 Proposition. Let $Lie(K) \to T(M)_x$ be the map $\xi\mapsto
\xi^\sharp_x$, and
$Lie(K)\to T(M)_x^*$ be its composition with the map $\imath_\omega$. Then the
transpose map is equal to
the differential $d\Phi_x:T(M)_x \to T(Lie(K)^*) = Lie(K)^*$
of $\Phi$ at $x$.

\proclaim  2.1.4 Corollary. The image of $d\Phi_x$ is equal to the annihilator
$Lie(K)_x^\perp$of the
Lie algebra $Lie(K)_x$ of the stabilizer $K_x$ of the point $x$.
 In particular, $x$ is not a critical point
of $\Phi$ if and only if $K_x$ is a finite group. The kernel of $d\Phi_x$ is
equal to the subspace
of all $v\in T(M)_x$ such that $\omega(v,\xi^{\sharp}_x) = 0$ for all $\xi \in
Lie(K)$.

\smallskip\noindent
2.1.5 We will be using the moment maps in the following situation. Let $G$ be a
reductive
algebraic group over the field of complex numbers ${\bf C}$ acting on a
projective nonsingular algebraic variety $X\subset {\bf P}(V)$ via a linear
representation on $V$.
We consider $G$ as a complex Lie group and consider the induced action of its
compact form $K$. One can choose a positive definite Hermitian inner product
$v\bullet w$ on $V$
such that $K$ acts on $V$ by means of a homomorphism
$\rho:K \to U(V)$. Now ${\bf P}(V)$ is equipped with the Fubini-Study
symplectic form, which is
 $U(V)$-invariant.

\smallskip\noindent
2.1.6. There is a canonical moment map ${\bf P}(V)\to {\rm Lie}(U(V))^*$. It
defines a canonical moment map $\Phi: M \to Lie(K)^*$.
 We want to describe its image
$\Phi (M)\subset Lie(K)^*$. Since $\Phi$ is $K$-equivariant, $\Phi (M)$
consists of
the union of co-adjoint orbits. We use some non-degenerate $K$-invariant
quadratic form on
$Lie(K)$ to identify $Lie(K)^*$ with $Lie(K)$ (in the case that $K$ is
semisimple, one simply
uses the Killing form).
This changes the moment map to a map $\Phi^*: M\to Lie(K)$.
Let $\hbar$ be a Cartan algebra of $Lie(K)$ and let $H$ be the maximal torus in
$K$ whose Lie algebra
is ${\rm Lie}(H) = \hbar$. Via derivations any character
$\chi: H \to U(1)$ of $H$ is identified with a linear function on $\hbar$. The
subset of such
functions is a lattice $\hbar_{\bf Z}$ in $\hbar$. We have $\hbar^* =
\hbar_{\bf Z}^*\otimes {\bf R}$. Under the isomorphism $\hbar \to
{\hbar}^*$ defined by the Killing form,
we obtain a ${\bf Q}$-vector subspace
$\hbar_{\bf Q}$ of $\hbar$ which defines a rational structure on $\hbar$. It
is known that each adjoint orbit intersects $\hbar$ at a unique orbit of the
Weyl group acting on $\hbar$. If we fix a positive Weyl
 chamber ${\hbar}_+$, then we obtain a bijective correspondence between adjoint
orbits and
 points of ${\hbar}_+$. This defines {\it the reduced moment map}:
$$
\Phi_{red}:M\to {\hbar}_+, x\mapsto \Phi^* (Kx)\bigcap {\hbar}_+.
$$

\smallskip\noindent
2.1.7 Let us go back to the algebro-geometric situation when $K$ is a compact
form
of a reductive group $G$ acting algebraically on a projective algebraic variety
$X \subset {\bf P}(V)$ via its linear representation
on $V$. A co-adjoint orbit ${\sl O}$ is called ${\it rational}$ if under the
correspondence
$$ \quad\hbox {co-adjoint orbits} \leftrightarrow \hbar_+$$
 it is defined by an element $\alpha$ in $\hbar_{\bf Q}$. We can write
$\alpha$ in the form ${\chi \over n}$ for some integer $n$ and $\chi \in {\cal
X}(T)$, where $T$ is a fixed maximal torus of $G$.
Elements of $\hbar_+$ of the form ${\chi \over n}$ will be called rational. We
denote the set of such elements by
$(\hbar_+)_{\bf Q}$. By the theorem of Borel-Weil,
${\chi}$ determines an irreducible representation $V(\chi)$ which can
realized as the space of sections
of an ample line bundle $L_\chi$ on the homogeneous space $G/B$, where $B$ is a
Borel
subgroup containing $T$. Its highest weight is equal to $\chi$.

\proclaim 2.1.8 Theorem. Let $\Phi_{red}:{\bf P}(V)\to {\hbar}_+$ be the
reduced moment map for the action of $K$
on ${\bf P}(V)$. Then $\Phi_{red}(X)$ is a compact convex subpolyhedron of
$\Phi_{red}({\bf P}(V))$ with vertices at rational points. Moreover
$$\Phi_{red}(X)\bigcap(\hbar_+)_{\bf Q} = \{{\chi \over n}:V(\chi)^*\quad
\hbox{ is a direct summand of $\Gamma (X,L^{\otimes n})$ as a $G$-module
$\}$}.$$

This result is due to D. Mumford (see the proof in {\bf [Ne2]}, Appendix. cf.
also {\bf [Br]}).
We refer to {\bf [MFK]}, Chapter 8, for the discussion of the similar convexity
result in the general situation of
2.1.6. Note that
the assertion is clear in the situation of 1.1.5. In fact, our calculation
shows that the
image of the moment map is equal to the convex hull $\overline {St(V)}$ of the
set
$$
St(V) = \{\chi: V_{\chi} \not= \{0\}\}.
$$

Since in our case $Lie(K)$ is abelian, we have $Lie(K) = \hbar_+$, so the
isomorphism $Lie(K) \cong Lie(K)^*$ identifies the image of the moment map with
the image of the reduced moment map. From this the assertion follows. This
convexity result was first observed by M. Atiyah {\bf
[At]}.

The following result from {\bf [Ne2]} gives the moment map interpretation of
the function $M^\bullet(x)$:

\proclaim 2.1.9 Theorem. For any $x \in X, M^L(x)$ is equal to the signed
distance from the origin to the boundary of $\Phi(\overline {G\cdot x})$.
If $x\in X^{us}(L)$, then the following three properties are equivalent:
\item {(i)} $M^L(x) = \|\Phi(x)\|$;
\item {(ii)} $x$ is a critical point for $|\!| \Phi |\!|^2$ with nonzero
critical value;
\item {(iii)} for all $t\in {\bf R},exp(t\Phi^*(x))\in G_x$, and the
complexification
of  $exp({\bf R}\Phi^*(x))$
is adapted for $M^L(x)$.

\medskip\noindent
{\bf 2.2} {\sl Relationship with geometric invariant theory}. We shall consider
the situation of 2.1.5.
Let $\Phi:X\to Lie(K)^*$ be the moment map as defined in 2.1.2. Let $L$ be the
restriction of
${\cal O}_{{\bf P}(V)}(1)$ to X. Recall that the stability of a point $x\in X$
with respect to $L$ is determined by the function $M^L(x)$ (see 1.1.7).
 Fix a norm $\|\ \|$ on $Lie(K)$ which is invariant with respect to the adjoint
action of $K$ on $Lie(K)$. The
following result follows from {\bf [KN]} but is explicitly stated in
L. Ness {\bf [Ne2]}, Theorem 2.2 (see also {\bf [Ki]}, 7.4,7.5):

\proclaim 2.2.1 Theorem.
\item {(i)} $ G \cdot \Phi^{-1}(0) = X^{ss}(L)_c:=\{x\in X^{ss}(L):G\cdot x$ is
closed $\}$;
\item {(ii)} $X^{ss}(L) = \{x\in X:\overline{G\cdot x}\cap \Phi^{-1}(0)\ne
\emptyset\}$;
\item{(iii)} $X^s(L) = X^{ss}(L)$ iff $0$ is not a critical value of $\Phi$.

\smallskip\noindent
2.2.2 By Theorem 2.2.1 (i), each closed orbit in $X^{ss}(L)$ has non-empty
intersection with the set $\Phi^{-1}(0)$. Since the latter is $K$-invariant,
this intersection is the union of $K$-orbits.
In fact, it consists of a unique orbit (see {\bf [KN]}, Theorem 0.1 (b)).
Since, set-theoretically, $X^{ss}(L)/\!/G$ is equal to the set of closed orbits
in
$X^{ss}(L)$, we obtain a homeomorphism
$$
\Phi^{-1}(0)/K \to X^{ss}(L)/\!/G.
$$

The orbit space $\Phi^{-1}(0)/K$ is called the {\it symplectic reduction} or
{\it Marsden-Weinstein reduction} of $X$ by $K$.
It has a natural structure of a symplectic manifold provided that
$0$ is not a critical value of $\Phi$ and $K$ acts freely on $\Phi^{-1}(0)$.
This is easily explained
with the help of the previous results. The first condition tells us that
$X^{ss}(L) = X^s(L)$, and the second one
implies that $X^{ss}(L)/\!/G = X^{s}(L)/G$ is a smooth projective variety. The
descent $M$ of $L$ to the quotient is an
ample line bundle. The cohomology class of the Marsden-Weinstein symplectic
form on $\Phi^{-1}(0)/K$ equals the first Chern class
of $M$. The Marsden-Weinstein construction gives a canonical way to descend a
representative of $c_1(L)$
to a representative of $c_1(M)$. This is done as follows. Let $x$ be a point in
$\Phi^{-1}(0)$ and
 $\bar x$ be the point of $\Phi^{-1}(0)/K$ corresponding to
the orbit $K\cdot x$. The tangent space at $\bar x$ is identified with the
quotient space
$Ker(d\Phi_x)/T(K\cdot x)_x$. By Corollary 2.1.4, $Ker(d\Phi_x)$ is orthogonal
to $T(K\cdot x)_x$ under $\omega_x$. Thus $\omega$ defines a skew-symmetric
bilinear form
on the quotient space. One checks that this defines a symplectic form on
$\Phi^{-1}(0)/K$.

\smallskip\noindent
2.2.3 Let $p\in \hbar_+$ be a point in the image of the moment map
$\Phi_{red}$, and
${\sl O}^p \subset Lie(K)^*$ be the corresponding co-adjoint orbit. It comes
equipped with a
canonical symplectic structure defined as follows. Given a co-adjoint orbit
${\sl O}\subset Lie(K)^*$
 and a point $q\in {\sl O}$, we identify the tangent space $T({\sl O})_q $ of
${\sl O}$ at $q$ with a subspace of $Lie(K)^*$. Now the skew-symmetric bilinear
map $(a,b)\mapsto q([a,b])$ is an element
$\omega_x$ of
${\Lambda}^2(T({\sl O})_q^*)$.

Let $\bar{\sl O}^p$ denote the symplectic manifold obtained from the symplectic
manifold
${\sl O}^p$ by replacing its symplectic form $\omega$ by $-\omega$. Then the
product
symplectic manifold $X\times \bar{\sl O}^p$ admits a moment map $ \Phi_p$
defined by the formula
$$
\Phi_p(x,q) = \Phi (x)-q.
$$
Now the set $\Phi_p^{-1}(0)$ becomes identified with  the set $\Phi^{-1}({\sl
O}^p)$ and
$$
\Phi_p^{-1}(0)/K = \Phi^{-1}(p)/K_p = \Phi^{-1}({\sl O}^p)/K
$$
where $K_p$ is the isotropy subgroup of $K$ in $K$ at $p$.
This quotient space is called the {\it Marsden-Weinstein reduction} or
{\it symplectic reduction} of $X$ with respect to $p$. Evidently it depends
only on the orbit of $p$.

The following result (see {\bf [Ne2]}, Appendix by D. Mumford) is another
interpretation of Theorem 2.1.8.

\proclaim 2.2.4 Theorem. Let $L$ be a very ample $G$-linearized line bundle on
$X$; $L({\chi \over n})$ denotes the line bundle on
$X\times G/B$ equal to the tensor product of the preimages
of the bundles $L^{\otimes n}$ and $L_\chi$ under the projection maps. Then
$L({\chi \over n})$ defines a $G$-equivariant embedding of $X\times G/B$ into
a projective space and  hence a moment map
$\Phi_{{\chi \over n}}:X\times G/B \to Lie(K)^*$. We have
$$\Phi_{{\chi \over n}}^{-1}(0) \cong \Phi^{-1}({\sl O}^{{\chi \over n}}).$$
In particular,
$$(X\times G/B)^{ss} (L({\chi\over n}))/\!/G = \Phi_{{\chi \over n}}^{-1}(0)/K
=
 \Phi^{-1}({\sl O}^{{\chi \over n}})/K.$$\par

\medskip\noindent
{\bf 2.3} {\sl Homological equivalence for $G$-linearized line bundles}. We
assume $X$ to be a projective variety.

\smallskip\noindent
2.3.1. Recall the definition of the Picard variety $Pic(X)_0$ and the
N\'eron-Severi group $NS(X)$.
We consider the Chern class map $c_1: Pic(X) \to H^2(X,{\bf Z})$ and put
$$Pic(X)_0 = Ker(c_1), NS(X) = Im(c_1).$$
One way to define $c_1$ is to choose a Hermitian metric on $L\in Pic(X)$, and
set $c_1(L)$ to be equal to the cohomology class
of the curvature form of this metric. Thus elements of $Pic(X)_0$ are
isomorphism classes of line bundles
which admit a Hermitian form with exact curvature form $\Theta$. If $\Theta$ is
given locally
as ${i \over 2\pi}d'd''log(\rho_U)$, then $\Theta$ is exact if and only if
there exists a global function
$\rho(x)$ such that $d'd''log(\rho_U/\rho)$ = 0 for each $U$. This implies that
$\rho_U/\rho = |\phi_U|$
for some holomorphic function
$\phi_U$ on $U$. Replacing the transition functions $\sigma_{UV}$ of $L$ by
$\phi_U \sigma_U \phi_V^{-1}$ we find an isomorphic bundle $L'$
with transition functions $\sigma'_{UV}$ with $|\sigma'_{UV}|$ = 1. Since
$\sigma'_{UV}$ are holomorphic, we obtain that the transition functions
of $L'$ are constants of absolute value 1. The Hermitian structure on $L'$ is
given by a global function
$\rho$.

\smallskip\noindent
2.3.2. Now assume $L$ is a  $G$-linearized line bundle, where $G$ is a complex
Lie group acting
holomorphically on $X$. Let
$K$ be its compact real form. By averaging over $X$ (here we use that $K$ and
$X$ are compact) we can find a Hermitian structure on $L$ such that
$K$ acts on $L$ preserving this structure (i.e. the maps $L_x\to L_{gx}$ are
unitary maps).
The curvature form of $L$ with respect to the $K$-invariant Hermitian metric is
a
$K$-invariant 2-form $\Theta$ of type (1,1). If $L$ is ample, $\Theta$ is a
K\"ahler form and its imaginary part is a
$K$-invariant symplectic form on $X$. A different choice of  $K$-invariant
Hermitian structure on $L$ replaces $\Theta$ with
$\Theta' = \Theta+{i\over 2\pi}d'd''\log(\rho)$ where $\rho$ is a $K$-invariant
positive-valued smooth real function
on $X$. This can be seen as follows. Obviously $\Theta-\Theta'$ is a
$K$-invariant 2-form $d'd''log\rho$ for some global function $\rho$. By the
invariance, we get for any
$g\in K, d'd''(log(\rho(gx)/\rho(x)) = 0$. This implies that $\rho(gx)/\rho(x)
= |c_g(x)|$ for some
 holomorphic function $c_g(x)$ on $X$. Since $X$ is compact, this function must
be constant, and the map
$g\mapsto |c_g|$ is a continuous homomorphism from $K$ to ${\bf R}^*$. Since
$K$ is compact and connected, it must be
trivial.
This gives $|c_g| = \rho(gx)/\rho(x) = 1$, hence $\rho(x)$ is $K$-invariant.

Taking the cohomology class of $\Theta$ we get a homomorphism:
$$c:Pic^G(X)\to H^2(X,{\bf Z}).$$

\proclaim 2.3.3 Proposition. The restriction of the canonical forgetful
homomorphism
$Pic^G(X)\to Pic(X)$ to $Ker(c)$ has its image equal to $Pic_0(X)$ and its
kernel isomorphic
to ${\cal X}(G)$. Moreover there is a natural section $s:Pic_0(X)\to Pic^G(X)$
defining an
isomorphism
$$  Ker(c) \cong {\cal X}(G) \times Pic(X)_0.$$

\proof
It is known (see {\bf [KKV]}) that the cokernel of the canonical forgetful
homomorphism
$Pic^G(X)\to Pic(X)$ is isomorphic to $Pic(G)$. Since the latter group is
finite and $Pic(X)_0$ is divisible, it is easy to see
that that $Ker(c)$ is mapped surjectively to
$Pic(X)_0$ with the kernel isomorphic to the group ${\cal X}(G)$ of characters
of $G$. To prove the assertion it suffices to construct a section
$s:Pic(X)_0\to Pic^G(X)$. Now let
$L\in Pic(X)_0$. We can find an open trivializing cover $\{U_i\}_{i\in I}$ of
$L$ such that  $L$ is defined by
constant transition functions $\sigma_{UV}$. For each
$g\in G$ the cover $\{g(U_i)\}_{i\in I}$ has the same property. Moreover,
$\sigma_{g(U)g(V)} = (g^{-1})^*(\sigma_{UV}) = \sigma_{UV}$. This shows that we
can define the ({\it trivial}) action of $G$ on the total space ${\bf L}$ of
$L$ by
the formula
$g\cdot (x,t)\to (g \cdot x, t)$
which is well-defined and is a $G$-linearization on $L$.
This allows one to define the section $s$ and check that
$s(Pic(X)_0)\cap {\cal X}(G) = \{1\}$.
\endproof

\proclaim 2.3.4 Definition. A $G$-linearized line bundle $L$ is called {\it
homologically trivial}
if $L\in Ker(c)$ and the image of $L$ in ${\cal X}(G)$ is the identity. In
other words, $L$ is
homologically trivial if it belongs to $Pic(X)_0$ as a line bundle, and the
$G$-linearization on $L$ is trivial (in the sense of the previous proof).
The subgroup of $Pic^G(X)$ formed by homologically trivial $G$-linearized line
bundles is denoted by $Pic^G(X)_0$. Two elements of $Pic^G(X)$ defining the
same element of the factor group
$Pic^G(X)/Pic^G(X)_0$ are called {\it homologically equivalent}.

\proclaim 2.3.5 Lemma. Let $L$ and $L'$ be two homologically equivalent
$G$-linearized line bundles. Then
for any $x\in X$ and any one-parameter subgroup $\lambda$
$$\mu^L(x,\lambda) = \mu^{L'}(x,\lambda).$$

\proof
 It is enough to verify that for any homologically trivial $G$-linearized line
bundle $L$ and any
one-parameter subgroup $\lambda$ we have $\mu^L(x,\lambda) = 1$. But this
follows immediately from the definition of the trivial $G$-linearization on
$L$.
\endproof

\proclaim 2.3.6 Theorem. Let $L$ and $L'$ be two ample $G$-linearized bundles.
Suppose that $L$ is homologically equivalent to
$L'$. Then the moment maps defined by the corresponding Fubini-Study symplectic
forms differ by a $K$-invariant constant.

\proof  Choose some symplectic curvature forms $\omega$ and $\omega'$ of $L$
and $L'$, respectively.
By assumption, $\omega ' = \omega + {i\over 2\pi}d'd''log(\rho)$ for some
positive-valued $K$-invariant function $\rho$ (see 2.3.2). Thus we can write
$\omega = \omega' +d{\theta}$,
where $\theta$ is a $K$-invariant 1-form of type $(0,1)$. By definition of the
moment map, we
have for each $\xi\in Lie(K)$
$$ d(\xi\circ\Phi') = \imath_{\omega'}(\xi^{\sharp}) = \imath_{\omega +
d\theta}(\xi^{\sharp}) =
d(\xi\circ\Phi)+L(\xi^\sharp)\theta-d\langle\xi^\sharp,\theta\rangle,$$
where $L(\xi^\sharp)$ denotes the Lie derivative  along the vector field
$\xi^\sharp$.
Since $\theta$ is $K$-invariant, we get $L(\xi^\sharp)\theta = 0$.
Because $G$ (hence also $K$) acts on $X$ holomorphically, each element of
Lie($G$) defines a holomorphic
vector field. In particular, $\xi^\sharp$ is holomorphic. But $\theta$ is of
type $(0,1)$, so $<\xi^\sharp, \theta> = 0$.
This shows that $ d(\xi\circ\Phi')  = d(\xi\circ\Phi)$ and hence
$\Phi' - \Phi$ is an  $ad(K)$-invariant constant in $Lie(K)^*$.
\endproof

\smallskip\noindent
2.3.7 Let $NS^G(X) = Pic^G(X)/Pic^G(X)_0$. By Proposition 2.3.3, we have an
exact sequence
$$0\to {\cal X}(G) \to NS^G(X) \to NS(X) \to A,$$
where $A$ is a finite group.
It is known that the N\'eron-Severi group $NS(X)$ is a finitely generated
abelian
group. Its rank is called the {\it Picard number} of $X$ and is denoted by
$\rho(X)$.  From this we infer that $NS^G(X)$ is a finitely
generated abelian group of rank $\rho^G(X)$ equal to
$\rho (X)+t(G)$, where $t(G)$ is the dimension of the radical $R(G)$ of $G$.
 Let
$$NS^G(X)_{\bf R} = NS^G(X)\otimes {\bf R}$$
be the finite-dimensional real vector space generated by $NS^G(X)$.

\smallskip\noindent
2.3.8 One can give the following interpretation of elements of $NS^G(X)_{\bf
R}$. First of all
the space $NS^G(X)_{\bf R}$ contains the lattice $Num^G(X)$ isomorphic to
$NS^G(X)/{\rm Torsion}$.
Secondly, it contains the ${\bf Q}$-vector space $NS^G(X)_{\bf Q}$. For each
$l\in NS^G(X)_{\bf Q}$ there is a smallest integer $n$
such that $nl$ is equal to the class $[L]$ of a $G$-linearized line bundle $L$
in
$Num^G(X)$. By using 2.3.2 and Theorem 2.3.6 we can interpret an element of
$NS^G(X)_{\bf R}$ as a moment map given by a symplectic forms which arises as
the imaginary part
of a closed real $K$-invariant 2-forms of type (1,1).

\medskip\noindent
{\bf 2.4} {\sl Stratification of the set of unstable points via moment map}.
Here we shall compare the Hesselink stratification
of the set $X\setminus X^{ss}(L)$ with the Ness-Kirwan stratification of the
same
 set using the Morse theory for the moment map.
The main references here are {\bf [Ki], [Ne2]}

\smallskip\noindent
2.4.1. Let $X$ be a compact symplectic manifold and $K$ be a compact Lie group
acting symplectically on it with  moment map
$\Phi: X \to Lie(K)^*$. We choose a $K$-invariant inner product on $Lie(K)$ and
identify $Lie(K)$ with $Lie(K)^*$.
Let $f(x) = \|\Phi(x)\|^2$.
It is obvious that the critical points of $\Phi$ are critical points
of  $f:X\to {\bf R}$. But some minimal critical points of $f$ need not  be
critical for
$\Phi$ (consider the case when 0 is a regular point of $\Phi$).
 For any $\beta\in Lie(K)$ we denote by $\Phi_{\beta}$ the composition
of  $\Phi$ and $\beta:Lie(K)^*\to {\bf R}$. Let $Z_{\beta}$ be the set of
critical points of $\Phi_{\beta}$ with critical value equal to
$\|\beta\|$. This is a symplectic submanifold of $X$ (possibly disconnected)
fixed by the subtorus $T_\beta$ equal to the closure in $H$ of the
real one-parameter subgroup $exp \; {\bf R}\beta$, where $H$ is a  maximal
subtorus of $K$.

The composition of $\Phi$ and the natural map
$Lie(K)^*\to Lie(H)^*$ is the moment map $\Phi_H$ for the action of $H$ on $X$.
If we use the
inner product on $Lie(K)$ to identify $Lie(H)$ with its dual, then $\Phi_H(x)$
equals the orthogonal projection of
$\Phi(x)$ to $Lie (H)$.  Thus, if $\Phi(x)\in Lie(H)$ then $\Phi_H(x) =
\Phi(x)$. Hence
such a point $x$ is a critical point for $f =\|\Phi\|^2$ if and only if it is a
critical point for the function
$f_H = \|\Phi_H\|^2$.  By a theorem of Atiyah {\bf [A]},
$\Phi_H (X)$ is equal to the convex hull of a finite set of points
$A$ in $Lie(H)$. Each point of $A$ is the image of a fixed point of $H$.
In general, the points of $A$ are
not rational points in the vector space $Lie(H)$.

\proclaim 2.4.2 Proposition.
\item{(i)} A point $x \in X$ is critical for $f$ if and only if the induced
tangent vector $\Phi^*(x)_x^\sharp$ at $x$
equals $ 0$.
\item{(ii)} Let $x \in X$ be a point such that $\Phi(x) \in \hbar$. Then $x$ is
a critical point for $f$
             if and only if it is a critical point for $f_H$.
\item{(iii)}
$x$ is a critical point for $f_H$ with critical value $\beta \ne 0$
if and only if $x\in Z_\beta$. In this case $\beta$ is the closest point to $0$
of the convex hull of some subset of the set $A$.

\proof This proposition follows from Lemmas 3.1, 3.3, and 3.12 from {\bf [Ki]}.
\endproof

\smallskip\noindent
2.4.3. Fix  a positive Weyl chamber $\hbar_+$ and
let ${\bf B}$ be the set of points $\beta$ in $\hbar_+$ which can be realized
as the closest point to the origin of the convex hull of some subset of the set
$A$. For each $\beta \in {\bf B}$ we set
$$C_\beta = K\cdot (Z_\beta \cap \Phi^{-1}(\beta)).$$
Let $S_\beta$ (resp. $Y_\beta$) denote the subset of points $x \in X$ such that
the limit set
of the path of steepest descent for $f$ from $x$
(with respect to some suitable $K$-invariant metric on $X$) is contained in
$C_\beta$ (resp. $Z_\beta$).
Then $Y_\beta$ is a locally closed submanifold of $S_\beta$.

\proclaim 2.4.4 Theorem. The critical set of $f$ is the disjoint union of the
closed subsets
${C_\beta,\beta\in {\bf B}}$.
The sets $\{S_\beta\}$ form a smooth stratification of $X$.

\proof This is Theorem 4.16 from {\bf [Ki]}.
\endproof

\smallskip\noindent
2.4.5 Now assume that the symplectic structure on $X$ is defined by the
imaginary part of
a K\"ahler form on $X$. Let $G$ be the complexification of $K$. This is a
reductive complex
algebraic group. We assume that the action of $K$ on $X$ is the restriction of
an action of $G$ on $X$ which
 preserves the K\"ahler structure of $X$. Then we have a stratification
$\{S_\beta\}$ chosen with respect
to the K\"ahler metric on $X$. In this case we can describe the stratum
$S_\beta$ as follows:
$$S_\beta = \{x\in X:\beta \  \hbox {is the unique closest point to 0 of}\
\Phi_{red} (\overline {G\cdot x})
\}.$$
Let $Z^{min}_\beta$ denote the minimal Morse stratum of $Z_\beta$ associated to
the function
$\|\Phi-\beta\|^2$ restricted to $Z_{\beta}$. Note that $\Phi_{\beta}$ is the
moment map for the symplectic manifold
$Z_\beta$ with respect to the stabilizer group $Stab_\beta$ of $\beta$ under
the adjoint action of $K$.
Let $Y^{min}_\beta$ be the pre-image of $Z^{min}_\beta$ under the natural map
$Y_{\beta} \to Z_\beta$.

For any
$\beta \in {\bf B}$ let
$$P(\beta)= \{g\in G:exp(it\beta)\cdot g\cdot exp(-it\beta) \hbox{ has a limit
in} \; G  \}.$$
It is a parabolic subgroup of $G$, and is the product of the Borel subgroup of
$G$ and $Stab_\beta$ which leaves $Y^{min}_{\beta}$ invariant.

\proclaim 2.4.6 Theorem.
\item{(i)} If $x \in Y^{min}_{\beta}$ then $\{ g \in G | gx \in
Y^{min}_{\beta} \} = P(\beta)$.
\item{(ii)} There is an isomorphism $S_{\beta} \cong
G\times_{P(\beta)}Y^{min}_{\beta}$;

\proof This is Theorem 6.18 and Lemma 6.15  from {\bf [Ki]}.
\endproof

\smallskip
\noindent
 The relationship between the stratification $\{S_\beta\}$ and the
stratification described in section 1.3 is as follows.
First of all we need to assume that $X$ is a projective algebraic variety and
the K\"ahler structure on $X$ is
given by the curvature form of an ample $L$ on $X$. In this case each $\beta$
is a rational vector in $\hbar_+ = Lie(H)_+$. For any such $\beta $ we can find
a positive integer
$n$ such that $n\beta$ is a primitive integral vector of $\hbar_+$ and hence
defines a primitive
one-parameter subgroup $\lambda_\beta$ of $G$.

\proclaim 2.4.7 Theorem. There is a bijective correspondence between the
moment map
strata $\{S_\beta\}_{\beta\not= 0}$ and the Hesselink strata
$\{S^L_{d,<\tau>}\}$ given by
$S_\beta \to S^L_{\|\beta\|,<\lambda_\beta>}$. Under this correspondence
\item {(i)} $P(\beta) = P(\lambda_\beta), Stab_\beta = L(\lambda_\beta)({\bf
R})$;
\item {(ii)}$Z_\beta = X_{\|\beta\|}^{\lambda_{\beta}}$; $Z^{min}_\beta =
Z^L_{\|\beta\|,\lambda_\beta}$;
\item {(iii)} $Y^{min}_\beta = S^L_{\|\beta\|,\lambda_\beta}$.

\proof See {\bf [Ki]}, \S12, or {\bf [Ne2]}.
\endproof

There is the following analog of Theorem 1.3.9:

\proclaim Theorem 2.4.8. \item {(i)} The set of locally closed subvarieties $S$
of $X$ which can be realized as the stratum
$S_\beta$ for some $K$-equivariant K\"ahler symplectic structure and $\beta \in
Lie(H)^*$ is finite.
\item {(ii)} The set of possible open subsets of $X$ which can be realized as a
stratum $S_0$ for some K\"ahler symplectic structure
on $X$ is finite.

\proof The second assertion obviously follows from the first one.  We prove the
first assertion by using induction on
the rank of $K$. The assertion is obvious if $K = \{1\}$. Applying induction to
the case when $X$ is equal to a connected component of $Z_\beta$ and $K
=Stab_\beta/T_\beta$, we obtain that the set
of open subsets of some connected component of $Z_\beta$ which can be realized
as the connected component of
the set $Z^{min}_\beta$ is finite. The finiteness of the set of subsets that
 can be realized as $Z_\beta$ follows from Lemma 1.3.6.
This implies that the set of closed subsets of $X$ which can be realized as the
subsets
$Z^{min}_\beta$ is finite. Now each $Z^{min}_\beta$ determines $Y^{min}_\beta$,
and  by Theorem 2.4.6,
 the stratum $S_\beta$.
\endproof

\medskip
\noindent
{\bf 2.5} {\sl K\"ahler quotients}.
We can extend many notions of GIT-quotients to the case of symplectic
reductions
with respect to a moment map $\Phi^\omega:X\to Lie(K)^*$
defined by a K\"ahler structure $\omega$ on
$X$ (see {\bf [Ki]}).

\smallskip\noindent
2.5.1 Perhaps the fastest way to go for a K\"ahler quotient is to define
$$X^{ss}(\omega)_c = G(\Phi^\omega)^{-1}(0). $$
The orbit space $X^{ss}(\omega)_c/G$ is Hausdorff.
In the case that $K$ acts freely on $(\Phi^\omega)^{-1}(0)$, the symplectic
reduction
$(\Phi^\omega)^{-1}(0)/K$ has a natural K\"ahler structure and
is diffeomeorphic to the orbit space $X^{ss}(\omega)_c/G$.

\smallskip\noindent
2.5.2 In general, for any $x\in X$ and $\lambda \in {\cal X}_*(G)$ we define
$$\mu^\omega (x,\lambda) = \|\lambda\| d_\lambda(0,\Phi^\omega (\overline
{G\cdot x})),$$
where $d_\lambda(0,A)$ denotes the signed distance
 from the origin to the boundary of the projection of
the set $A$ to the positive ray spanned by $\lambda$.
Thus we can define
$$M^\omega(x) = \sup_\lambda d_\lambda(0,\Phi^\omega (\overline {G\cdot x}))
$$
so that $M^\omega(x)$ is equal to the signed
distance from the origin to $\Phi^\omega(\overline {G\cdot x})$,
and set
$$X^{ss}(\omega) = \{x\in X:M^\omega(x) \leq 0 \},$$
$$X^s(\omega) = \{x\in X:M^\omega(x) < 0 \},$$
and
$$X^{sss}(\omega) = X^{ss}(\omega)\setminus X^s(\omega).$$
One can show that $X^{ss}(\omega)_c$ is the union of closed orbits in
$X^{ss}(\omega)$
or equivalently $X^{ss}(\omega)$ is the union of those orbits in X
whose closures meet $X^{ss}(\omega)_c$.
 As in the case of Hodge K\"ahler structures
we can define the set of primitive adapted one-parameter
subgroups $\Lambda^\omega(x)$. There are analogs of Theorems 1.2.4 and 1.2.6 in
our situation.

\smallskip\noindent
2.5.3 The analog of the categorical quotient $X^{ss}/\!/G$ is the orbit space
$(\Phi^\omega)^{-1}(0)/K$. The analog of the geometric quotient is the orbit
space
$X^s(\omega)\cap (\Phi^\omega)^{-1}(0)/K$.

\bigskip
\noindent
{\bf \S 3. The G-ample cone}.
\smallskip\noindent
{\bf 3.1} {\sl $G$-effective line bundles}. We will assume that
$X$ is projective and normal.

\proclaim 3.1.1 Definition. A $G$-linearized line bundle $L$ on $X$ is called
{\it $G$-effective} if
$X^{ss}(L) \not= \emptyset$. A $G$-effective ample $G$-linearized line bundle
is called {\it $G$-ample}.

\proclaim 3.1.2 Proposition. Let $L$ be an ample $G$-linearized line bundle.
The following assertions are equivalent:
\item{(i)} $L$ is $G$-effective;
\item{(ii)}$\Gamma (X,L^{\otimes n})^G \not= \{0\}$ for some $n > 0$;
\item{(iii)} If $\Phi: X \to Lie(K)^*$ is the moment map associated to a
$G$-equivariant embedding of $X$
into a projective space given by some positive tensor power of $L$, then
$\Phi^{-1}(0) \not= \emptyset$.

\proof
(i)$\Leftrightarrow $(ii). Follows from the definition of semi-stable points.

\noindent
(ii) $\Leftrightarrow $ (iii). Follows from Theorem 2.2.1.
\endproof

\proclaim 3.1.3 Proposition. If $L$ and $M$ are two $G$-effective
$G$-linearized bundles then $L\otimes M$ is $G$-effective.

\proof
Let $x\in X^{ss}(L)\bigcap X^{ss}(M)$; then there exists $s\in \Gamma
(X,L^{\otimes n})^G$ and
$s'\in \Gamma (X,M^{\otimes m})^G$, for some $n,m > 0$ such that $s(x) \not= 0,
s'(x)\not= 0$
and
$X_s$ and $X_{s'}$ are both affine. Without loss of generality we may assume
that $n = m$. Then the $G$-invariant section
$s\otimes s'$ of $L^{\otimes n}\otimes M^{\otimes n}$
does not vanish at $x$. Using the Segre embedding we get that $X_{s\otimes s'}$
is affine. This shows that
$x\in X^{ss}(L\otimes M)$, hence $L\otimes M$ is $G$-effective.

\endproof

\proclaim 3.1.4 Proposition. Let $L,L'$ be two ample $G$-linearized line
bundles. Suppose they are homologically equivalent.
Then $X^{ss}(L)$ = $X^{ss}(L')$.

\proof This follows immediately from the numerical criterion of stability and
Lemma 2.3.5.
\endproof

\smallskip\noindent
3.1.5 {\sl Remark}. One should compare this result with Corollary 1.20 from
{\bf [MFK]}. Under the assumption that $Hom(G,{\bf C}^*) = \{1\}$ it asserts
that $X^s(L) = X^s(L')$ for any $G$-linearized
line bundles defining the same element in $NS(X) = Pic(X)/Pic(X)_0.$

\smallskip\noindent
3.1.6 An element of $NS^G(X)$ will be called {\it $G$-ample} if it can be
represented by
an $G$-effective ample line bundle. By Proposition 3.1.4,
all ample $G$-linearized line bundles in the same
$G$-ample homological equivalence class are $G$-effective.

Let $NS^G(X)^+$ denote the subset of $G$-ample homological equivalence classes.
Using Proposition 3.1.3, one checks that  it is a semigroup in $NS^G(X)$.

\smallskip\noindent
3.1.7 Using the terminology from 2.5.2, we can extend the previous results and
the definitions to
any point $[\omega]$ from $NS^G(X)_{\bf R}$ (not necessary rational). By
definition $[\omega]$ consists of a cohomology class of some closed
$K$-invariant 2-form $\omega$ of type $(1,1)$ with non-degenerate imaginary
part ${\rm Im}(\omega)$ and a choice of a moment map defined by the
symplectic form ${\rm Im}(\omega)$. The proof of Theorem 2.3.6 and the
discussion of 2.3.8
 show that the set $X^{ss}([\omega])$ is well defined.
It equals the union of orbits which contain in its closure a point from
$G\cdot(\Phi^\omega)^{-1}(0)$, where $\omega$ is a representative of
$[\omega]$.
We say that $[\omega]$ is $G$-effective if $X^{ss}([\omega])\ne \emptyset$. By
the convexity of a moment map the set of $G$-effective points is closed under
addition. It is obviously closed under multiplication by positive scalars. Thus
we can define the
convex cone $EF^G(X)$ of $G$-effective points in $NS^G(X)_{\bf R}$. We shall
call it the {\it $G$-effective cone} of $X$.

\medskip
\noindent
{\bf 3.2} {\sl The G-ample cone}. Here comes our main definition:

\proclaim 3.2.1 Definition. The {\it $G$-ample cone} (for the action of $G$ on
$X$) is the convex
cone in $NS^G(X)_{\bf R}$ spanned by the subset $NS^G(X)^+$. It is denoted by
$C^G(X)$.

\smallskip\noindent
3.2.2 Let $X$ be a compact complex K\"ahler manifold. The subset of the space
$H^{1,1}(X,{\bf R})$ formed by the classes of K\"ahler forms is an open convex
cone. It is called the
{\it K\"ahler cone}. Its integral points are the classes of Hodge K\"ahler
forms.
By a theorem of Kodaira, each such class is the first Chern class of an ample
line bundle $L$.
The subcone of the K\"ahler cone spanned by its integral points is called the
{\it ample cone} and is denoted by
$A^1(X)^+$. It is not empty if and only if  $X$ is a projective algebraic
variety.
It spans the subspace $A^1(X)^+_{\bf R}$ of $H^{1,1}(X,{\bf R})$
 formed by the cohomology classes of algebraic cycles of codimension 1. The
dimension
of this subspace is equal to the Picard number of $X$.
The closure of the ample cone consists of the classes of numerically effective
(nef) line
bundles {\bf
[Kl]}.
Recall that a line bundle is called {\it nef} if its restriction to any curve
is an
effective line bundle. Under the forgetful map $NS^G(X)_{\bf R}\to NS(X)_{\bf
R}$, the
$G$-ample cone is mapped
to the closure of the ample cone. Of course they are equal if $G = \{1\}$.
Summing up we conclude that
$$C^G(X) = EF^G(X)\cap \alpha^{-1}(A^1(X)^+_{\bf R}),$$
where $\alpha:NS^G(X)_{\bf R}\to NS(X)_{\bf R} \to H^{1,1}(X,{\bf R})$ is the
composition map.
\smallskip\noindent
3.2.3 {\sl Remark}. There could be no $G$-effective
 $G$-linearized line bundles, so $C^G(X)$ could be empty. The simplest example
is any homogeneous space $X = G/P$, where $P$ is a parabolic subgroup of a
reductive group $G$.

\smallskip\noindent
3.2.4 Let us recall some terminology from the theory of convex sets (see {\bf
[Ro]}). A {\it convex set} is a subset $S$ of a real vector space $V$
such that for any $x,y\in S$ and any nonnegative real numbers $\lambda,\mu$
with
$\lambda + \mu = 1$, $\lambda x + \mu y \in S$. A convex set is a {\it convex
cone} if additionally
for each $x\in S$ and any nonnegative real number $\lambda$, $\lambda x \in S$.
The minimal dimension of an affine subspace containing a convex cone $S$ is
called the dimension of $S$ and is denoted by $dim(S)$.
The {\it relative interior} $ri(S)$ of $S$ is the interior of $S$ in the sense
of the topology of the minimal affine space containing
$S$. The complementary set $r\partial (S) := S\ \setminus ri(S)$ is called the
{\it relative boundary} of $S$.
A closed convex cone is called {\it polyhedral} if it is equal to the
intersection of a finite number of closed half-spaces.
An affine hyperplane $H$ is called a {\it supporting hyperplane} at a point
$x\in \partial (S)$ if $x\in S$ and $S$ lies in one of the two half-spaces
defined by this hyperplane. This half-space is called the {\it supporting
half-space} at
$x$. A supporting hyperplane of a convex cone is clearly a linear hyperplane. A
{\it face} of a closed convex set $S$ is the intersection of $S$ with a
supporting hyperplane. Each face is a closed convex set. A function
$f:V \to {\bf R}\cup \{\infty \}$ is called
{\it lower convex} if $f(x+y) \leq f(x)+f(y)$ for any $x,y \in V$. It is called
{\it positively homogeneous} if
$f(\lambda x) = \lambda f(x)$ for any nonnegative $\lambda$. If $f$ is a lower
convex positively homogeneous function,
then the set $\{x\in V : f(x) \leq 0\}$ is a closed convex cone. For any set
$A$ we denote by $\langle A \rangle$ the smallest closed convex cone containing
$A$.

\proclaim 3.2.5 Lemma.
\item {(i)} The relative interior of any non-empty convex set of positive
dimension is nonempty;
\item {(ii)} A point belongs to the boundary of a convex set $S$ if and only if
there exists a supporting hyperplane at this point;
\item{(iii)} Each closed convex set is equal to the intersection of the
supporting half-spaces chosen for each point of its boundary;
\item {(iv)} The closure of a convex set is a convex set.

\smallskip
Now we can go back to our  convex set $C^G(X)$.

\proclaim 3.2.6 Lemma. For each $x\in X$ the function $Pic^G(X) \to {\bf Z},L
\to M^L(x)$ factors through $NS^G(X)$ and can be uniquely
extended to a positively homogeneous lower convex function
$M^\bullet(x):NS^G(X)_{{\bf R}} \to {\bf R}$.

\proof By Lemma 2.3.5, $M^L(x) = M^{L'}(x)$ if $L$ is homologically equivalent
to $L'$. This shows that we can descend
the function $L\to M^L(x)$ to the factor group $NS^G(X)$. By using 1.1.1 (ii),
we find
$$M^{L^{\otimes n}}(x) = nM^L(x)$$
for any nonnegative integer $n$, and
$$M^{L\otimes L'}(x) = sup_\lambda {\mu^{L\otimes L'}(x,\lambda) \over
\|\lambda\|}
=sup_\lambda ({\mu^L(x,\lambda) + \mu^{L'}(x,\lambda) \over \|\lambda\|}) $$
$$\leq sup_\lambda {\mu^L(x,\lambda) \over \|\lambda\|}+sup_\lambda
{\mu^{L'}(x,\lambda) \over \|\lambda\|}=
M^L(x)+M^{L'}(x).$$
Let us denote by [$L$] the class of $L\in Pic^G(X)$ in $NS^G(X)$. If we put
$$M^{{1\over n}[L]}(x) = {1\over n}M^L(x) \quad\hbox {if $n$ is positive},$$
$$ M^{{1\over n}[L]}(x) = {1\over -n}M^{L^{-1}}(x) \quad\hbox {if $n$ is
negative},$$
then we will be able to extend the function $L\to M^L(x)$ to a unique function
on
$$NS^G(X)_{\bf Q} =NS^G(X)\otimes {\bf Q}$$ satisfying
$$ M^{l+l'}(x) \leq M^l(x)+M^{l'}(x),$$
$$M^{\alpha l}(x) = {\alpha} M^l(x)$$
for any $l,l'\in NS^G(X)_{\bf Q}$ and any non-negative rational number
$\alpha$.
Now let us choose a basis ${e_1,\ldots,e_n}$ in $NS^G(X)_{\bf Q}$ and define a
norm  by
$\|\sum_ia_ie_i\| = sup_i|a_i|$. Then for any $l = \sum_ia_ie_i$ we have
$$M^l(x) \leq \sum_iM^{a_ie_i}(x)\leq
\sum_isup\{|a_i|M^{e_i}(x),|a_i|M^{-e_i}(x)\} \leq C\|l\|$$
for some constant $C$ independent of $l$. This implies that
$$M^{l \over \|l\|} (x) \le C.$$
But $M^{l \over \|l\|} (x)$ is clearly bounded from below by the minimum of
the continuous function ${l \over \|l\|}: \rightarrow {\mu^{l \over \|l\|}(x,
\lambda) \over \|\lambda\|}$
defined over a sphere of radius 1 for any fixed $\lambda$.
Since a bounded convex function is continuous in $NS^G(X)_{\bf Q}$,
the above argument allows us to extend the function $M^\bullet (x)$ to a unique
continuous
function on $NS^G(X)_{\bf R}$.
Obviously it must be lower convex and positively homogeneous.
\endproof

\proclaim 3.2.7 Proposition. Assume ${\rm dim} \; C^G(X) = {\rm dim} \;
NS^G(X)_{\bf R}$. For each $x\in X$ the restriction of the function $M^\bullet
(x)$ to $C^G(X)$ coincides
with the function $[\omega] \to M^\omega(x)$, where $\omega$ is any
representative of $[\omega]$ and
$M^\omega(x)$ is defined in 2.5.2.

\proof
The two functions coincide on the set of rational points in $C^G(X)$.
\endproof

\smallskip
What can we say about the boundary of $C^G(X)$? It is clear that any integral
point in it is
 represented either by an ample line $G$-bundle
or by a nef but non-ample
$G$-bundle $L$.

\proclaim 3.2.8 Proposition. Let $L$ be an ample $G$-linearized line bundle in
the boundary of $C^G(X)$. Then
$L\in C^G(X)$ and $X^s(L) = \emptyset$.

\proof
First let us show that $X^s(L) = \emptyset$. Suppose $x\in X^s(L)$, then
$M^L(x) < 0$ and hence the intersection of the open set
$\{[\omega]:M^{[\omega]}(x) <0\}$ with the open cone
$\alpha^{-1}(A^1(X)^+_{\bf R})$ is an open neighborhood of $[L]$ contained in
$C^G(X)$. This contradicts the assumption that $L$ is on the boundary.

Now let us show that $L\in C^G(X)$, i.e., $X^{ss}(L)\ne \emptyset.$ Suppose
$X^{ss}(L) =\emptyset$. This means that $M^L(x) > 0$ for all $x\in X$.  By
Theorem 1.3.9
the set ${\cal S}$ of open subsets
$U$ of $X$ which can be realized as the set $X^{ss}(M)$ for some ample
$G$-linearized line bundle $M$ is finite.
Choose a point $x_U$ from each such $U$. Then $X^{ss}(L) = \emptyset$ if and
only if
$M^L(x_U) > 0$ for each $U$. This shows that $L$ is contained in the
intersection of
the open set $\{[\omega]:M^{[\omega]}(x_U) > 0, U\in {\cal S}\}$ with the open
cone
$\alpha^{-1}(A^1(X)^+_{\bf R})$. Thus $L$ belongs to the complement of the
closure of
$C^G(X)$, contradicting the assumption on $L$.
\endproof

\proclaim 3.2.9 Corollary. Assume that there exists an ample $G$-linearized
line
bundle $L$ with
$X^s(L)\not=\emptyset$. Then the interior of $C^G(X)$ is not empty.

\medskip\noindent
3.3 {\sl  Walls and chambers}.

\proclaim 3.3.1 Definition. A subset $H$ of $C^G(X)$ is called a {\it wall}
if there exists a point $x\in X_{(>0)}:=\{x:{\rm dim}G_x > 0\}$
such that $H = H(x):= \{l\in C^G(X):M^l(x) = 0\}$. A wall is called {\it
proper} if it
is not equal to the whole cone $C^G(X)$. A connected component of the
complement of the union of walls in $C^G(X)$ is called
a {\it chamber}.

\smallskip

\proclaim 3.3.2 Theorem. Let $l,l'$ be two points from $C^G(X)$.
\item {(i)}  $l$ belongs to some wall if and only if $X^{sss}(l) \not=
\emptyset $;
\item {(ii)} $l$ and $l'$ belong to the same chamber if and only if $X^s(l) =
X^{ss}(l) =  X^{ss}(l') = X^s(l')$;
\item {(iii)} each chamber is a convex cone.

\proof (i) If $l$ belongs to some wall then there exists a point $x\in X$ such
that $M^l(x) = 0$.
By 2.5.2, $x\in X^{sss}(l)$. Conversely,
if $x\in X^{sss}(l)$ then the closure of $G\cdot x$ contains a closed orbit
$G\cdot y$
in $X^{sss}(l)$ with stabilizer
$G_y$ of positive dimension. Thus $l$ lies on the wall $H(y)$.

\noindent (ii) If $l$ belongs to a chamber $C$, then $X^s(l) = \{x\in X: M^L(l)
< 0\}$.
The function $M^\bullet (x)$ does not change sign in the interior of $C$.
Because if it does, we can find a point $l'$
in $C$ such that $M^{l'}(x) = 0$. Arguing as in the proof of (i) we find that
$l'$ belongs to a wall.
This shows that the set $X^s(l)$
does not depend on $l$. We shall denote it by $X^s(C)$.  Let us prove the
converse. Assume $X^s(C) = X^s(C')$ for two chambers
$C$ and $C'$.
We have $C, C'\subset \cap_{x\in X^s(C)}\{M^\bullet(x) < 0\}$. Pick  $l$ in
$\cap_{x\in X^s(C)}\{M^\bullet(x) < 0\}$.
Then $X^s(C) \subset X^s(l)$, hence we get $X^s(C)/G \subset X^s(l)/G$. By (i),
the first quotient is compact. Therefore
the second quotient is compact and $X^s(l) = X^{ss}(l)$. This shows that $l$
does not belong to the union
of walls. Since $\cap_{x\in X^s(C)}\{M^\bullet(x) < 0\}$ is convex, and hence
connected, it must coincide with  $C$.
This shows that $C = C'$ and the set
$X^s(C)$ determines the chamber $C$ uniquely.

\noindent (iii) This obviously follows from the fact that the function
$M^\bullet(x)$ does not take the value zero in a chamber and is positively
homogeneous lower convex.
\endproof

\proclaim 3.3.3 Theorem. There are only finitely many walls.

\proof
 For any wall $H$ let $X(H) = \{x\in X_{(>0)}:H = H(x)\}$ where $X_{(>0)}$ is
the set of points in $X$ with positive dimensional
stabilizer. Obviously
$$ X(H) = \cap_{l\in H} (X^{ss}(l))\cap X_{(>0)}).$$
By Theorem 2.4.8, we can find a finite set of points $l_1,\ldots,l_N$ in
$C^G(X)$
such that
for any $l\in C^G(X)$, the set $X^{ss}(l)$ equals one of the sets
$X^{ss}(l_i)$.
Using the above description of $X(H)$ we obtain that there are only finitely
many subsets of $X$ which are equal to the subset
$X(H)$ for some wall $H$. However, two walls $H,H'$ are equal if and only if
$X(H) = X(H')$ if and only if $X(H)\cap X(H') \not= \emptyset$.
This proves the assertion.
\endproof

\proclaim 3.3.4 Proposition. Each proper wall is
equal to a countable union of convex closed cones of dimension strictly less
than the dimension of $C^G(X)$.

\proof
 A wall is a set of the form
$\{ l \in C^G(X)_{\bf R}: M^l(x) = 0 \}$. The set  $\{ l \in C^G(X)_{\bf R}:
M^l(x) \le 0 \}$ is a
convex closed cone.
The function $M^\bullet(x)$ is defined as the supremum of linear functions
$\mu^l(x,\lambda)/||\lambda||$. Those of them which are not identically zero
and {\it vanish} at $l$
 are obviously among the supporting hyperplanes of the cone.
It follows from the properties of these
functions stated in 1.1.1 that there are only countably many different
functions among them
(see {\bf [MFK]},p.59). It remains to use that the boundary of a convex closed
cone $C$
is equal to the  union of
intersections of $C$ with its supporting hyperplanes.
\endproof

\smallskip
Since a wall is a countable union of convex cones we can  speak about its
dimension
(being the maximum of the dimensions of
the cones).

\proclaim 3.3.5 Proposition. Assume all walls are proper. Then
$X^s(l)\ne\emptyset$ for any $l$ in the interior of
$C^G(X)$.

\proof
Assume $l$ is contained in an open subset $U$
of $C^G(X)$. By Theorem 3.3.3,
the union of walls is a proper closed subset of $C^G(X)$. If $l$ does not
belong to any
wall, the assertion is obvious. If it does, we choose a line through $l$ which
contains two
points $l_1,l_2 \in U$ which do not lie in the union of walls. Then
$X^s(l_1)\cap X^s(l_2)\ne \emptyset$ so that
$M^{l_1}(x) < 0$ and $M^{l_2}(x) < 0$ for some $x\in X$. By convexity, $M^l(x)
< 0$, hence $x\in X^s(l)$.
\endproof

\smallskip\noindent
Together with Proposition 3.2.8 this gives the following

\proclaim 3.3.6 Corollary. Assume all walls are proper.
The subset of the boundary of $C^G(X)$
whose image under the forgetful map is contained in the ample cone
equals the set
$\{l\in C^G(X): X^s(l) =\emptyset\}.$

\proclaim 3.3.7 Lemma. Let $H$ be a wall. There exists a point $x\in X$ such
that $H = H(x)$ but
$H \not= H(y)$ for any $y\in \overline {G\cdot x}\setminus G\cdot x$.

\proof Take any point $x\in X$ such that $H=H(x)$. Let $y\in \overline {G\cdot
x}\setminus G\cdot x$.
Assume $H = H(y)$. Since ${\rm dim} \;G\cdot y < {\rm dim}\; G\cdot x$, we can
replace
$x$ by $y$ and continue this process until we find the needed point.
\endproof

\proclaim 3.3.8 Definition. A point $x\in X$ is called a {\it pivotal point}
for a wall $H$ if
$H = H(x)$ and $H\not= H(y)$ for any $y\in \overline {G\cdot x}\setminus G\cdot
x$.

\proclaim 3.3.9 Lemma. Let $\Phi: X \to Lie(K)^*$ be a moment map associated to
some symplectic structure
on $X$. For any $x\in \Phi^{-1}(0)$ the connected component of $G_x$ is a
reductive algebraic group.
In particular, if $G\cdot x$ is a closed orbit in $X^{ss}(l)$ for some $l \in
C^G(X)$, then
the connected component of $G_x$ is a reductive algebraic group.

\proof This is Lemma 2.5 from {\bf [Ki]}. In the case when the moment map is
defined
by an ample line bundle $L$ the proof is as follows. By Theorem 2.2.1, $G\cdot
x$ is closed in $X^{ss}(L)$.
 Let $\pi:X^{ss}(L)\to X^{ss}(L)/\!/G$ be the quotient projection. Its fibers
are affine, and hence
$G\cdot x$ is a closed subset of an affine variety. Thus $G\cdot x = G/G_x$ is
an affine variety.
 Now the assertion follows from
{\bf [B-B1]}.
\endproof

\proclaim 3.3.10 Proposition. Let $x$ be a pivotal point for a wall $H$. Then
the connected component of $G_x$ is a reductive subgroup of $G$.

\proof Let $l\in H\cap C^G(X)$ but $l\notin H(y)$ for any $y\in \overline
{G\cdot x}\setminus G\cdot x$.
Since the number of walls is finite, we can always find such an $l$.
Now, by the choice of $l$,
one sees that $G\cdot x$ is a closed orbit in $X^{ss}(l)$
 and we can apply the previous lemma to obtain that
the connected component of $G_x$ is reductive.
\endproof

\proclaim 3.3.11 Definition.
Let ${\cal S}$ be the set of all possible moment map stratifications of $X$
 defined by points $l\in C^G(X)$. There is a natural partial order in ${\cal
S}$  by refinement.
We say that $s$ refines $s'$ if any stratum of $s'$ is a union of strata of
$s$.
In this case, we sometimes write $s < s'$ if $s \ne s'$.

Let $C$ be a subset of $C^G(X)$.
For any $s \in {\cal S}$, let $C_s$ be the set of elements $l\in C$ which
define the stratification $s$.
Then $C = \bigcup_{s\in{\cal S}}C_s$.

The next proposition implies that $s$ refines any element in the closure of
$C_s$.

\proclaim 3.3.12  Proposition. Let $l_n \in C^G(X)$ be a sequence of points in
$C^G(X)$
that induce the same stratification $s \in {\cal S}$. Assume that $l_n \to l
\in C^G(X)$.
Then $s$ refines the  stratification $s'$ induced by $l$.

\proof  Let $S_{\beta (l_n)}$ be a stratum in $s$.
Since $S_{\beta (l_n)}= GY^{min}_{\beta (l_n)}$,
it suffices to show that $Y^{min}_{\beta (l_n)}$ is contained in a stratum of
$s'$.
Let $\lambda(l_n)$ be the one-parameter subgroup generated by $\beta (l_n)$.
Then, by Theorem 2.4.7, we can also write $Y^{min}_{\beta (l_n)}$ as $S_{d_n,
\lambda(l_n)}$,
where $d_n = |\!|\beta (l_n)|\!|$.
By passing to a subsequence, we may assume that
$d_n \to d$,
$\beta (l_n) \to \beta (l)$.
Let $\lambda(l)$ be the one-parameter subgroup generated by  $\beta (l)$.
For any point $x \in Y^{min}_{\beta (l_n)}$,
we have $M^{l_n}(x) = \bar{\mu}^{l_n}(x, \lambda(l_n))= d_n$ (see Theorem
2.1.9 (iii)). By taking limits, we find that
$M^{l}(x) = \bar{\mu}^l(x, \lambda(l)) = d$.
If $d=0$, we obtain that $Y^{min}_{\beta (l_n)} \subset X^{ss}(l)$.
If $d\ne0$, we obtain that $Y^{min}_{\beta (l_n)} \subset S_{d, <\lambda(l)>}$.
This completes the proof.
\endproof

\smallskip
The following will play the key role in the sequel.

\proclaim 3.3.13 Lemma. Let $l_0$ be in the closure of a chamber $C$. Then
there is a sequence of points $l_n$ in $C$ which induce
the same stratification $s$ and $l_n \to l_0$.
In particular, $X^{ss}(l_0)$ is a union of the strata of $s$.

\proof   Take a basis $U_1 \supset U_2 \supset \cdots \supset U_n \supset
\cdots$
of open subsets containing $l_0$.
We have that $U_n \cap C \ne \emptyset$ and is open in $C$. Write ${\cal S}
=\{s_1, \cdots, s_m \}$.
If $s_1$ occurs as the induced stratification for some point $l_n$ in $U_n \cap
C$ for every  $n$,
then we are done. Suppose it does not. Then there is $n_1$ such that $(U_{n_1}
\cap C)_{s_1} = \emptyset$.
Consider $s_2$. If $s_2$ occurs as the induced stratification for some point
$l_n$ in $U_n \cap C$ for every
$n \ge n_1$, then we are done. Suppose it does not. Then there is $n_2 \ge n_1$
such that $(U_{n_2} \cap C)_{s_2}=\emptyset $. The next step is to consider
$s_3$, and so on.
Since the set ${\cal S} =\{s_1, \cdots s_m \}$ is finite and
$U_1 \cap C \supset U_2 \cap C \supset \cdots \supset U_n \cap C \supset
\cdots$ is infinite,
there must be an $s_i \in {\cal S}$ such that $s_i$ occurs as the induced
stratification
for some point $l_n$ in $U_n \cap C$ for every $n \ge n_{i-1}$.

Then $\{l_n\}$ ($n \ge n_{i-1}$) is the desired sequence.
\endproof

\smallskip\noindent
3.3.14
For any $x\in X$ and $L\in Pic^G(X)$ let $\rho_x(L):G_x\to GL({\bf L}_x) \cong
{\bf C}^*$ be
the isotropy representation. For any $\lambda\in {\cal X}_*(G_x)$ the
composition
$$\langle \lambda,\rho_x(L)\rangle = \rho_x (L) \circ \lambda : {\bf C}^* \to
{\bf C}^*$$
is equal to the map $t \to t^{\mu^L(x,\lambda)}$.
The correspondence $L\to \rho_x(L)$ defines a homomorphism
$$\rho_x:Pic^G(X) \to {\cal X}(G_x).$$
Obviously $\rho_x(L)$ is trivial if $L$ is $G$-homologically trivial. Thus we
can define a homomorphism
$$\bar\rho_x:NS^G(X)\to {\cal X}(G_x)$$
to be the one naturally induced by  $\rho_x:Pic^G(X) \to {\cal X}(G_x).$
By linearity we can extend $\bar\rho_x$ to a linear map
$${\bar\rho}_x^{\bf R}:NS^G(X)_{\bf R}\to {\cal X}(G_x)\otimes {\bf R}.$$
For any $l\in NS^G(X)_{\bf R}$ and $\lambda\in {\cal X}_*(G_x)$ we have, by
continuity,
$$\langle\lambda,\bar\rho_x^{\bf R}(l)\rangle (t) = t^{\mu^l(x,\lambda)}.$$

We denote by ${\bf E}_x$ the kernel of the map ${\bar\rho}_x^{\bf R}$.

\proclaim 3.3.15 Proposition.
Let $x$ be a pivotal point of a wall $H$. Assume that all walls are proper.
Then  ${\bf E}_x \not= NS^G(X)_{\bf R}$.

\proof Let $l$ be chosen as in the proof of Proposition 3.3.10. Since all walls
are proper,
we can find an open neighborhood $U$ of $l$ in $NS^G(X)_{\bf R}$
such that it contains a point $l'\in C^G(X)$ which is not contained in any wall
and the stratification induced by $l'$ refines
 the stratification induced by $l$ (see Lemma 3.3.13).
 We can also find a continuous path in $U$ which starts at
$l'$, ends at $l$ and does not cross any walls. Since $M^\bullet(x)$ is a
continuous function, we
have $M^{l'}(y) < 0$ implies $M^l(y) \leq 0$, and $M^l(y) < 0$ implies
$M^{l'}(y) < 0$. This shows that
$$X^s(l') \subset X^{ss}(l),\   X^s(l) \subset X^s(l').$$

In particular, $X^{us}(l) \subset X^{us}(l')$. Let us now use the moment map
stratification
for $X$ with respect to  $l$ and $l'$. We can write
$$X = X^{ss}(l)\cup X^{us}(l) = X^{s}(l')\cup S_{\beta_1}\cup \ldots
\cup S_{\beta_k}\cup X^{us}(l).$$
Here,
we use the fact that the strata of $X^{us}(l)$ coincide with the union of some
strata of $X$ with respect to
$l'$. Now the point $x$ must belong to some stratum $S_{\beta_i}$
since it is semi-stable with respect to
$l$ but unstable with respect to $l'$. Let us denote $\beta_i$ by $\lambda$.
Let
$y = \lim_{t\to 0}\lambda (t)\cdot x$. It is contained in $S_\lambda \subset
X^{ss}(l)$
and is fixed by  $\lambda ({\bf C}^*)$ (see 1.3.2). By the choice of $l$ it
must belong
to the orbit of $x$. Now by Theorem 1.2.6, $\lambda \in \Lambda^{l'}(y)$ hence
$M^{l'}(y) = {\mu^{l'}(y,\lambda)\over \|\lambda\|} > 0$.
 By Theorem 2.4.6 (i) $y = g\cdot x$ for some $g\in P(\lambda)$.
Applying the known properties of the function $\mu^{l'}(x,\lambda)$ (see
1.1.1), we obtain
$$\mu^{l'}(y,\lambda)  = \mu^{l'}(g \cdot x, \lambda) =
\mu^{l'}( x,g^{-1}\circ \lambda\circ g) = \mu^{l'}( x, \lambda )> 0 ,$$
 hence $l'\notin {\bf E}_x$.
\endproof

\proclaim 3.3.16 Corollary. Assume all walls are proper. Let $H$ be a wall and
let
$x\in H$ be its pivotal point. Then
\item {(i)} $G_x$ is a reductive group with non-trivial radical $R(G_x)$;
\item{(ii)} $H  \subset {\bf E}_x;$
\item {(iii)} if ${\rm codim}H = 1$ then $H$ equals the closure of its set of
rational points.

\proof The first assertion obviously follows from the fact that ${\bf E}_x
\not= NS^G(X)$ implies
${\cal X}(G_x) \not= \{1\}$.

Let us prove (ii). Take any $l \in H(x)$. Since $x \in X^{ss}(l)$,
for any one-parameter subgroup $\lambda$ in $G_x$ we have
$\mu^l(x,\lambda) = 0$. In fact, otherwise either
$\mu^l(x,\lambda)$ or $\mu^l(x,\lambda^{-1})$ is positive,
 and then $x\in X^{us}(l)$.
This shows that the composition of $\lambda \colon {\bf C}^*\to G_x$
 with $\bar\rho_x^{\bf R}(l)\colon G_x \to{\bf C}^* $ is trivial. Hence $l\in
{\bf E}_x$.

To prove (iii) we use that, by (ii), $H$ is contained in ${\bf E}_x$. Since
${\bf E}_x$ is a proper subspace, $H$ spans it. This implies that the relative
interior $ri(H)$ of $H$ contains an open subset
of ${\bf E}_x$. Since ${\bf E}_x$ can be defined by a linear equation over
${\bf Q}$ this open subset contains a
dense subset of rational points.
\endproof

\proclaim 3.3.17 Proposition.
Assume that all walls are proper. Then any wall is contained in a codimension 1
wall.

\proof Since the union of walls is a closed subset of $C^G(X)$, each chamber is
an open subset
of $C^G(X)$. Suppose there is a wall $H$ of codimension $\geq 2$ which is not
contained in any codimension 1 wall.
Let $l$ be a point in the relative interior of $H$. Then there exists an open
subset $U$ of $C^G(X)$
containing $l$ such that $U\setminus H\cap U$ is contained in the complement of
the union of walls. Since
$H$ is of codimension $\geq 2$, the set $U\setminus H\cap U$ is connected.
Hence there exists a chamber
$C$ containing this set. Let $l_1,l_2 \in C$ such that
$l$ is on the line joining $l_1$ with $l_2$. As in the proof of Proposition
3.3.5,
$$X^s(l_1)\cap X^s(l_2)  \subset  X^s(l).$$
By Theorem 3.3.2 (ii), we see that the left-hand side is equal to
$X^{ss}(l_1)$.
Since the quotient $X^{ss}(l_1)//G = X^s(l_1)/G$ is compact, $X^s(l)/G$ is
compact, and hence
$X^{ss}(l) = X^s(l)$. But this contradicts
the assumption that $l$ belongs to a wall.
\endproof

\proclaim 3.3.18 Theorem. Assume that all walls are proper.
Then the intersection of the boundary of a chamber with interior of
$C^G(X)$ is contained in a union of finitely many rational supporting
hyperplanes.

\proof
Let $l$ be an interior point of $C^G(X)$ which belongs to the boundary of a
chamber $C$. Let $H_1,\ldots ,H_N$ be codimension 1 walls containing $l$. We
claim
that, for at least one of these walls $H = H_i$, the closure $\bar C$ of $C$
contains a
point in the relative interior of $H_i$.
In fact,
otherwise $l$ is contained in an open neighborhood $U$ such that the
intersection of $U$
with the union of chambers is connected (see the proof of the previous
Proposition).
Since $U$ is not contained in $C$ this contradicts the definition of chambers.
  By Proposition 3.3.15 and Corollary 3.3.16, $H$ spans a rational hyperplane
${\bf E}_x$, where
$x$ is a pivotal point of $H$. This hyperplane is the needed supporting
hyperplane for the closure $\bar C$ of $C$ at the point
$l$. To verify this it suffices to show that $C\cap {\bf E}_x = \emptyset$.
Suppose $l'$ is in the intersection.
Choose a point $l''$ in the relative interior of $H$ which belongs to the
boundary of $C$. Then the segment of the line joining the points
$l'$ with $l''$ is contained in ${\bf E}_x$ which is spanned by $H$. On the
other hand it contains only one point from
$H$, the end point $l''$. Since  $l''\in ri(H)$,
this is obviously impossible.
\endproof

\smallskip\noindent
3.3.19 { \sl Remark}. Theorem 3.3.18
 shows that each chamber is a rational convex polyhedron away from the boundary
of $C^G(X)$. There is no reason to expect that the boundary of $C^G(X)$ is
polyhedral. This is not true even when
$G$ is trivial. However, some assumptions on $X$, for example, $X$ is a Fano
variety (the dual of the canonical line bundle is ample)
may imply that the closure of $C^G(X)$ is a convex polyhedral cone.

\smallskip\noindent
3.3.20 { \sl Example}.
 Let $G$ be an $n$-dimensional torus $({\bf C}^*)^n$, acting on $X = {\bf
P}(V)$ via  a linear representation
$\rho:G \to GL(V)$.
Then $NS^G(X) \cong Pic(X)\times {\cal X}(G)\cong {\bf Z}^{n+1}$. The splitting
is achieved by fixing the $G$-linearization on
${\cal O}_X(1)$ defined by the linear representation $\rho$.
%Since $X^{ss}(L)$ does not change when we raise $L$ into a positive tensor
%%power,
In addition, $NS^G(X)^+ = {\bf Z}_{>0}\times NS^G(X)_1$, where $NS^G(X)_1$ is
the group of $G$-linearizations on
the line bundle ${\cal O}_X(1)$ identified with ${\cal X}(G)$. Thus $L_0$ =
$({\cal O}_X(1),\rho)$ corresponds to the zero in
${\cal X}(G)$. By 1.1.5, $X^{ss}(L_0) \not= \emptyset $ (resp.$X^{s}(L_0) \not=
\emptyset$) if and only if $0 \in St(V)$ (resp.
$0\in int(St(V))$ (see  2.1.8 for the notation $St(V)$). If $L_\chi$ is defined
by twisting the linearization of
$L_0$ by a character $\chi$, we obtain that $X^{ss}(L_\chi) \not= \emptyset $
(resp.$X^{s}(L_\chi) \not= \emptyset$) if and only
if $\chi^{-1} \in St(V)$ (resp.$\chi^{-1}\in int(St(V))$). In particular, we
get that $C^G(X)$ is equal to the cone
over the convex hull $\overline {St(V)}$ of $St(V)$. Comparing this with 2.1.8,
we obtain that $C^G(X)$ is equal to the cone over
the image of the moment map for the Fubini-Study symplectic form on $X$.

For any $x\in X$ the image of the orbit closure
$\overline{G\cdot x}$ under the moment map is equal to the
convex hull $P(x)$ of the state set
$st(x)$ of $x$. This is a subpolytope of $\overline {St(V)}$ of codimension
equal to $\dim G_x > 0$. A wall $H(x)$ is
the cone over $P(x)$ with $\dim G_x > 0$. The union of walls is equal to the
cone over the set of critical points of
the moment map.

\smallskip\noindent
3.3.21 { \sl Example}. Let $G= SL(n+1)$ act diagonally on $X = ({\bf P}^n)^m.$
We have
$Pic(X)\cong Pic^G(X)\cong {\bf Z}^m$. Each line bundle $L$ over $X$ is
isomorphic to the line bundle
$L_{\bf k} = p_1^*({\cal O}_{{\bf P}^n}(k_1))\otimes\ldots \otimes p_m^*({\cal
O}_{{\bf P}^n}(k_m))$
for some
 ${\bf k} = (k_1,\ldots,k_m)\in {\bf Z}^m$. Here
$p_i:X\to {\bf P}^n$ denotes the projection map to the $i$-th factor.
It is easy to see that $L_{\bf k}$ is nef (resp. ample) if and only if all
$k_i$ are nonnegative (resp. positive) integers.
 Let ${\cal P} = (p_1,\ldots,p_m)\in X$. Using the numerical criterion of
stability, one easily verifies
that
${\cal P}\in X^{ss}(L_{\bf k})$ if and only if for any proper linear
subspace $W$ of ${\bf P}^n$
$$ \sum_{i,p_i\in W}k_i \le {\dim W+1 \over n+1}(\sum\limits_{i=1}^mk_i).$$
The strict inequality characterizes stable points. This easily implies that
$$X^{ss}(L_{\bf k}))\ne \emptyset \Leftrightarrow (n+1){\rm max}_i\{k_i\}\le
\sum\limits_{i=1}^mk_i.$$
Let $\Delta_{n,m} = \{x = (x_1,\ldots,x_m)\in {\bf R}^m:\sum\limits_{i=1}^mx_i
= n+1, 0\le x_i\le 1, i = 1,\ldots, m\}.$
This is the so-called $(m-1)$-dimensional {\it hypersimplex} of type $n$.

Consider the cone over $\Delta_{n,m}$ in ${\bf R}^{m+1}$
$$C\Delta_{n,m} = \{(x,\lambda)\in {\bf R}^m\times {\bf R}_+:x\in \lambda
\Delta_{n,m}\}.$$
We have the injective map
$$ Pic^G(X) \to {\bf R}^{m+1}, L_{\bf k}\mapsto
(k_1,\ldots,k_m,(n+1)^{-1}\sum\limits_{i=1}^mk_i),$$
which allows us to identify $Pic^G(X)$ with a subset of ${\bf R}^{m+1}$. We
have
$$Pic^G(X)\cap C\Delta_{n,m} = \{\ L\in Pic^G(X): L \ \ \hbox{is
nef},X^{ss}(L)\ne \emptyset\}.$$
Thus the $G$-ample cone $C^G(X)$ coincides with the interior of
$C\Delta_{n,m}$.

Observe that ${\cal P}\in X^{sss}(L_{\bf k})$ if and only if there exists a
subspace
$W$ of dimension $d, 0\le d \le n-1$ such that
$$(n+1)\sum_{i,p_i\in W}k_i = (\dim W +1)\sum\limits_{i=1}^mk_i.$$
This is equivalent to the condition that $L_{\bf k}$ belongs to the hyperplane
$$H_{I,d}:=\{(x_1,\ldots,x_m,\lambda)\in {\bf R}^m:\sum_{i\in I}x_i = \lambda
d\},$$
where $I$ is a non-empty subset of $\{1,\ldots,m\}$.
Thus a chamber is a connected component of $C\Delta_{n,m}\setminus
\bigcup_{I,d}H_{I,d}$. A
wall is defined by intersection of the interior of $C\Delta_{n,m}$ with a
subspace
$$H_{I_1,\ldots,I_s,d_1,\ldots,d_s}:=\{(x_1,\ldots,x_m,\lambda)\in {\bf
R}^m:\sum_{i\in I_j}x_i = \lambda d_j, j=1,\ldots,s\},$$
where $\{1,\ldots,m\} = I_1\coprod\ldots\coprod I_s, d_1+\ldots+d_s = n+1-s,
s\ge 2.$

Note that the set $\Delta_{n,m}$ is the image of the moment map for the natural
action of the
torus $T= ({\bf C}^*)^m$ on
${\bf P}(\bigwedge^{n+1} {\bf C}^m)$. Comparing with Example 3.3.20, we obtain
that
the closure of
$C^{SL(n+1)}(({\bf P}^n)^m)$ is equal to $C^T({\bf P}(\bigwedge^{n+1} {\bf
C}^m))$. There
is a reason for this. For any ${\cal P} = (p_1,\cdots,p_m)\in ({\bf P}^n)^m$
one can
consider the matrix
$A$ of size $(n+1)\times m$ whose $i$-th column is a vector in ${\bf C}^{n+1}$
representing
the point $p_i$. Let $E(A)$ be the point of the Grassmann variety
$G(n+1,m)\subset {\bf P}(\bigwedge^{n+1} {\bf C}^m)$ defined by the matrix $A$.
A different choice of
coordinates of the points $p_i$ replaces $E(A)$ by the point $t\cdot E(A)$ for
some $t\in T$.
%Replacing ${\cal P}$ by $g\cdot {\cal P}, g\in SL(n+1)$ does not change the
%%point $E(A)$.
In this way we obtain a bijection between $SL(n+1)$-orbits of points
$(p_1,\ldots,p_m)\in ({\bf P}^n)^m$ with $\langle p_1,\ldots,p_m\rangle = {\bf
P}^n$  and
orbits
of $T$ on $G(n+1,m)$. This is called the {\it Gelfand-MacPherson
correspondence}.
The $T$-ample cone $C^T(G(n+1,m))$ equals $C^T({\bf P}(\bigwedge^{n+1} {\bf
C}^m)$ and
hence coincides with the closure of the $SL(n+1)$-ample cone of $({\bf
P}^n)^m$.
The Gelfand-MacPherson correspondence defines a natural isomorphism between the
two
GIT-quotients corresponding to the same
point in $C\Delta_{n,m}$.

\eject
\noindent
{\bf 3.4. GIT-equivalence classes}.

\proclaim 3.4.1 Definition. Two  elements $l$ and $l'$ in $C^G(X)$ are called
GIT-equivalent
(resp. weakly GIT-equivalent) if
$X^{ss}(l) = X^{ss}(l')$ (resp. $X^s(l) = X^s(l')$).

\smallskip
Let $E \subset C^G(X)$ be a GIT-equivalence class.
We denote by $X^{ss}(E)$ the subset of $X$ equal to $X^{ss}(l)$ for any
$l\in E$. Clearly for any $l\in E$ the subset $X^s(l)$ is equal to the subset
of those points in
$X^{ss}(E)$ whose orbits are closed in $X^{ss}(E)$ and whose stabilizers are
finite.
This shows that this set is independent of a choice of $l$, so we can denote it
by $X^s(E)$.
 In particular, GIT-equivalence implies weak GIT-equivalence.
 Similarly we introduce the subset $X^{ss}(E)_{(>0)}$ of points in $X^{ss}(E)$
whose stabilizers are of positive dimension. It is clear that
$$E\cap H(x) \ne \emptyset \Leftrightarrow E\subset H(x) \Leftrightarrow x\in
X^{ss}(E)_{(>0)}$$
where $x$ is a pivotal point for the wall $H(x)$.
Examples of GIT-equivalence classes are chambers (Theorem 3.3.2). For these
equivalence classes $X^{ss}(E) = X^s(E)$.
Chambers are the only GIT-equivalence classes which coincide with weak
GIT-equivalence classes.

\proclaim 3.4.2 Theorem.
Let $H$ be a wall and $H^\circ$ be the subset of points from $H$ which do not
lie in any
other wall unless it contains $H$ entirely. Assume $H^\circ\ne \emptyset$.
Then each connected component of $H^\circ$ is contained in a $GIT$-equivalence
class.
Conversely, for each
equivalence class $E$  which is not a chamber there exists a wall $H$ such that
$E$ is the union of some connected components of $H^\circ$.  In particular, a
subset is
a  connected component of  a $GIT$-equivalence class if and only if
it is a chamber or a  connected component of  $H^\circ$ for some wall $H$.

\proof
 Obviously, $l,l'\in H^\circ$ if and only if $X^{sss}(l) = X^{sss}(l').$
For any $x\in X$ the function $M^\bullet(x)$ either vanishes identically on
$H^\circ$ or
does not take the value zero at any point.
This implies that this function
does not change its sign at any connected component $F$ of $H^\circ$. Thus for
any $l,l'\in F$,
$X^s(l) = X^s(l')$. Together with $X^{sss}(l) = X^{sss}(l')$ we obtain that $F$
is contained
in a GIT-equivalence class.
By the previous discussion, any equivalence class $E$ not equal to a chamber
is contained in $H^\circ$ for some $H$. Hence it must be equal to
a union of connected components of $H^\circ$.
\endproof

\proclaim 3.4.3 Lemma. Let $l$ and $l'$ be two points in $C^G(X)$.
If $X^{ss}(l) \subset X^{ss}(l')$, then $X^{s}(l') \subset X^{s}(l)$.
Consequently, if $X^{ss}(l) \subset X^{ss}(l')$ and $X^{s}(l) \subset
X^{s}(l')$,
then $l$ and $l'$ are weakly equivalent.

\proof
If $X^{s}(l') = \emptyset$, then there is nothing to prove.
Assume that $X^{s}(l') \ne \emptyset$. If $x\in X^s(l')$ then $G\cdot x$ is
closed
in $X^{ss}(l')$ and hence
in $X^{ss}(l)$. In addition $G_x$ is finite, so $x\in X^s(l)$. This proves the
inclusion
$X^s(l')\subset X^s(l).$
\endproof

\proclaim 3.4.4 Theorem. Let $E$ be a GIT-equivalent class which is contained
in some wall $H$. Then its convex hull ${\rm CH}(E)$ is contained in a weak
GIT-equivalent class.
Moreover ${\rm CH}(E)\cap H^\circ \subset E$.

\proof
Let $l,l'\in E$ be joined by a segment $S$ of a straight line
and $l''$ be a point on this segment.
By the lower convexity of the functions $M^\bullet(x)$ we have
$X^{s}(l)\cap X^{s}(l') \subset X^s (l'')$
and $X^{ss}(l)\cap X^{ss}(l') \subset X^{ss}(l'')$.
Since
$X^{ss}(l) = X^{ss}(l')$ and hence $X^{s}(l) = X^{s}(l')$,
 we obtain $X^{s}(l) \subset X^{s}(l'')$,
 $X^{ss}(l) \subset X^{ss}(l'')$.
By the previous lemma, $l$ and $l''$
are weakly GIT-equivalent.
Assume that $l''$ does not lie on any other wall unless it contains $H$
entirely,
 i.e., $l''\in H^\circ$.
 Then the proof of the previous theorem shows that $X^{ss}(l) = X^{ss}(l'')$,
hence $l,l''\in E$.
This shows that $S\cap H^\circ \subset E$. Hence ${\rm CH}(E) \cap E \subset
H^\circ$.
\endproof

\proclaim 3.4.5 Definition. A non-empty subset of $C^G(X)$ is called a {\it
cell} if it is a chamber or
a connected component of $H^\circ$ for some wall $H$.

By Theorem 3.4.2, one can give an equivalent definition by defining a cell to
be a connected component of
a GIT equivalence class.

\proclaim 3.4.6 Definition. A pair of chambers
$(C, C')$ are called {\it relevant} to a cell $F$ if $F \subset \overline{C}
\cap \overline{C'}$ and
there is a {\it straight} path $l: [-1, 1] \rightarrow C^G(X)$ such that
$l([-1, 0)) \subset C$, $l(0) \in F$ and $l((0, 1]) \subset C'$.

\proclaim 3.4.7 Proposition.
 Let $(C, C')$ be a pair of chambers  relevant to a cell $F$. Then,
$$X^s(F) = X^s(C) \cap X^s(C'),  X^{ss}(F) \supset X^s(C) \cup X^s (C'). $$

\proof
 The fact that $X^{ss}(F) \supset X^s(C) \cup X^s (C')$ and
$X^s(F) \supset X^s(C) \cap X^s(C')$ follows from the lower convexity of the
function
$M^\bullet(x)$. The fact that $X^s(F) \subset X^s(C) \cap X^s(C')$
follows from the continuity of $M^\bullet(x)$.
\endproof

\medskip
\noindent
{\bf 3.5} {\sl Abundant actions}.
\smallskip\noindent
\proclaim 3.5.1 Definition. We say that $Pic^G(X)$ (or the action) is {\it
abundant} if for any pivotal point
$x$ of any wall the isotropy homomorphism
$\rho_x:Pic^G(X)\to {\cal X}(G_x)$ has finite cokernel.

\smallskip
The abundance helps to control the codimensions of walls.

\proclaim 3.5.2 Proposition. Assume that $Pic^G(X)$ is abundant. Then for any
pivotal point
$x$ of a codimension k wall $H$ the radical of the
stabilizer group $G_x$ is of dimension k or less.

\proof This follows immediately from Proposition 3.3.16.
\endproof

\proclaim 3.5.3 Theorem.  Let $B$ be a Borel subgroup of $G$ and let $G$ act on
$X \times G/B$ diagonally.
Then $Pic^G(X \times G/B)$ is abundant. In particular,
when $G$ is a torus, then $Pic^G(X)$ is abundant.

\proof
Notice first that $Pic^G(G/B) \cong {\cal X}(T)$, where $T$ is a maximal torus
contained in
$B$ (see {\bf [KKV]}, p.65).
We want to show that for any point $(x, g[B])$ with reductive stabilizer
the isotropy representation homomorphism
$$\rho_{(x, g[B])}: Pic^G(X \times G/B) \rightarrow {\cal X}(G_{(x, g[B])})$$
is surjective.
Obviously $G_{(x, g[B])} = G_x \cap G_{g[B]}$.
By conjugation, we may assume that
$R_{(x, g[B])} \subset R_x \subset T$, where $R_x$ denotes the radical of the
stabilizer subgroup $G_x$. Pick any character $\chi \in {\cal X}(G_{(x,
g[B])})$.
Extend it to a character $\tilde{\chi} \in {\cal X}(T)$.
Let $L_{\tilde\chi}\in Pic^G(G/B)$
be the line bundle asscoated to the character $\tilde \chi$,
 and let $L$ be a point in its inverse image in $Pic^G(X \times G/B)$ under the
projection  map $X\times G/B\to G/B$.
Then  $\rho_{(x, g[B])}(L) = \tilde\chi|_{R_{(x, g[B])}} = \chi|_{R_{(x,
g[B])}}.$
Thus $\rho_{(x, g[B])}$ is surjective.
Hence $Pic^G(X \times G/B)$ is abundant.

The last statement follows immediately because $G=B$ when $G$ is a torus.
\endproof

\bigskip\noindent
{\bf  \S 4. Variation of Quotients.}

\smallskip\noindent
{\bf 4.1} {\sl Faithful walls}.

\proclaim 4.1.1 Definition. A codimension 1 wall $H$ is called {\it faithful}
if for any pivotal point $x$ with $H(x) = H$ the radical of $G_x$ is
one-dimensional,
{\it truly faithful} if  the stabilizer $G_x$ is a one-dimensional torus.
%{\it perfect} if it is truly faithful and in addition all such stabilizers
%are conjugate to each other.
A cell $E$ is called faithful
(truly faithful) if it is contained in a faithful (truly faithful) wall.

%\smallskip\noindent
%4.1.2 Example  Let $G$ be a reductive algebraic group of rank 1 (e.g.,
%%$G=SL(2)$) acting on a
%projective variety $X$. Assume that all walls are proper. Then
%each wall is perfect. In fact by Proposition 3.3.15 the stabilizer
%of its pivotal point has non-trivial radical.
%Since $G$
%is of rank 1, all proper reductive subgroups must be tori,
% all tori are one-dimensional and  conjugate to each other.

\proclaim 4.1.2 Proposition.  Let $G$ be a reductive algebraic group acting on
a
nonsingular projective variety $X$. Assume that $Pic^G(X)$ is abundant.
Then all codimension 1 walls are  faithful.

\proof Follows from Corollary 3.3.16.
\endproof

\proclaim 4.1.3 Corollary. Let $G$ act on $X \times G/B$. Then
all codimension 1 walls are truly faithful.

\proof  Let $(x, g[B]) \in X \times G/B$. Then
$G_{(x, g[B])} = G_x \cap G_{g[B]} \subset G_{g[B]} = gBg^{-1}$.
Now if $(x, g[B])$ is a pivotal point for a  codimension 1 wall,
then $G_{(x, g[B])}$ is reductive.
This implies that $G_{(x, g[B])}$ is a torus.
Now the corollary follows from Theorem 3.5.3.
\endproof

\smallskip
Recall that GIT quotients of $X \times G/B$ by $G$
can be identified with symplectic reductions
of $X$ by $K$. This corollary will assure that our main theorem on variation of
quotients
applies to symplectic reductions for general coadjoint orbits --- and this is
exactly what
motivates this paper.

\smallskip
It worths mentioning that if we take $G$ to be a torus, then we get

\proclaim 4.1.4 Corollary. Let $G$ be a torus. Then
all codimension 1 walls are truly faithful.

\smallskip As remarked in the introduction, the new feature in this corollary
is that
it takes into account the variation of moment maps as well as the characters of
the
torus.

\medskip
\noindent
{\bf 4.2  Variation of quotients.}
\smallskip
In this section we assume that $X$ is nonsingular and all walls in $C^G(X)$ are
proper.

\proclaim 4.2.1 Lemma.  Let $F$ be a cell contained in
the closure of  another cell   $E$.
Then
\item{(i)}  $X^{ss}(E) \subset X^{ss} (F)$;
\item{(ii)}  $X^s(F) \subset X^s(E)$;
\item{(iii)}  the inclusion $X^{ss}(E) \subset X^{ss} (F)$ induces a morphism
$X^{ss}(E)/\!/G  \rightarrow X^{ss} (F)/\!/G$ which is an isomorphism over
$X^s(F)/G$.

\proof  (i) Let $l\in E$ and $x\in X^{ss}(l) = X^{ss}(E)$. Then $M^l(x)\le 0$
so by continuity
$M^{l'}(x)\le 0$ for any $l'\in F$. This proves that $X^{ss}(E)\subset
X^{ss}(l') = X^{ss}(F)$.

\noindent
(ii) In the previous notation we have $M^{l'}(x) < 0$ for any $x\in X^{s}(F)$.
Thus, by continuity,
$M^l(x) < 0$ for any $l\in E$ and hence $x\in X^s(E)$.

\noindent
(iii) Follows from (i) and (ii).
\endproof
\smallskip
\noindent
4.2.2 Let $(C^+,C^-$) is a pair of chambers relevant to
a cell $F$. Let $l_0 $ be a point in $F$.
Lemma 3.3.13 implies  that we can choose $l^+ \in C^+$
and $l^- \in C^-$ such that their induced
stratifications (cf. 1.3 and 2.5) of $X$
can be arranged  as follows:
$$X = X^{ss}(l_0) \cup X^{us}(l_0),$$
$$X = X^s(l^+) \cup S^{l^+}_{\alpha_1} \cup \cdots \cup S^{l^+}_{\alpha_p} \cup
X^{us}(l_0),$$
and
$$X = X^s (l^-) \cup S^{l^-}_{\beta_1} \cup \cdots \cup  S^{l^-}_{\beta_q} \cup
X^{us}(l_0).$$
To simplify the notation we shall assume that each $Z^{min}_{\alpha_i}$ or
$Z^{min}_{\beta_j}$
is connected.

Let
$$f^+:X^s(C^+)/G\to X^{ss}(F)/\!/G,\ \ f^-:X^s(C^-)/G\to X^{ss}(F)/\!/G$$
be the morphisms defined in Lemma 4.2.1. They are birational morphisms of
projective varieties
which are isomorphisms over the subset $X^s(F)/G$ of $X^{ss}(F)/\!/G$. The goal
of this
section is to describe the fibres of the morphisms $f^+$ and $f^-$.

\proclaim 4.2.3 Lemma.
Keep the previous notation and assume that the cell $F$ is truly faithful;
then we have that all $G$-orbits of points from  $Z^{min}_{\alpha_i}$
are closed in $X^{ss}(l_0)$ and all non-stable
closed orbits  in $X^{ss}(l_0)$ meet some $Z^{min}_{\alpha_i}$.
In addition, up to conjugation, $\alpha_i$ form the set of one-parameter
subgroups (without parametrizations) of $G$ that have nonempty
fixed point set on $X^{sss}(l_0)_c$, where
$X^{sss}(l_0)_c$ is the union of closed orbits in $X^{sss}(l_0)$.
Similar statements are also true
for $\beta_j$ and $Z^{min}_{\beta_j}$.

\proof
Since $F$ is truly faithful, $X^{ss}(l_0)$ does not contain points with
stabilizer of
dimension $> 1$. This implies that all $G$-orbits of points from
$Z^{min}_{\alpha_i}$
are closed in $X^{ss}(l_0)$.
Now, by 4.2.2,
$$X^{ss}(l_0) = X^s(l^+) \cup S^{l^+}_{\alpha_1} \cup \cdots \cup
S^{l^+}_{\alpha_p}.$$
The points of $X^s(l^+)$ have finite isotropy subgroups, so they do not lie in
$X^{ss}(l_0)_c \setminus X^s(l_0)$. Clearly points from $S^{l^+}_{\alpha_i}
\setminus
Z^{min}_{\alpha_i}$ are not in closed orbits of $X^{ss}(l_0)$.
Hence all non-stable
closed orbits  in $X^{ss}(l_0)$ must meet some $Z^{min}_{\alpha_i}$. This shows
that
$$X^{ss}(l_0)_c = X^s(l_0) \cup \bigcup_i (G \cdot Z^{min}_{\alpha_i}). $$
Since $F$ is a truly faithful cell, the isotropy subgroup of any point from
$Z^{min}_{\alpha_i}$
is exactly the one-parameter subgroup generated by $\alpha_i$.
Thus the  last statement for $\alpha_i$ follows readily from the above
observation.
The assertions for  $\beta_j$ can be proved similarly.
\endproof

\proclaim 4.2.4 Lemma.
Keep the previous notation and assume that the cell $F$ is truly faithful; then
we have
\item{(i)} $p = q$;
\item{(ii)} Under a suitable arrangement (including choosing suitable Weyl
chambers),
$$ Z^{min}_{\alpha_i}= Z^{min}_{\beta_i},\ 1\le i \le q;$$
\item{(iii)} $\alpha_i = - c_i\beta_i$ ($\ 1\le i \le q$) for some positive
number $c_i$;
 \item{(iv)} $S^{l^-}_{\beta_j} \cap S^{l^+}_{\alpha_i} = \emptyset$ if $j \ne
i$.

\proof  (i) and (ii) follow from the proof of Lemma 4.2.3 because
$$X^{ss}(l_0)_c = X^s(l_0) \cup \bigcup_i (G \cdot Z^{min}_{\alpha_i}) $$
and
$$X^{ss}(l_0)_c = X^s(l_0) \cup \bigcup_i (G \cdot Z^{min}_{\beta_i}). $$

Now let us prove (iii).
We have already seen that $\alpha_i$ and $\beta_i$ generate the same subgroup
of $G$, i.e., they differ by a constant multiple. Now by the constructions of
the stata, we  have
$S^{l^+}_{\alpha_i} = G \cdot Y^{min}_{\alpha_i}$  and $S^{l^-}_{\beta_i} = G
\cdot Y^{min}_{\beta_i}$,
where $Y^{min}_{\alpha_i}$ is the preimage over $Z^{min}_{\alpha_i}$ under the
Bialynicki-Birula
contraction determined by $\alpha_i$. Similarly $Y^{min}_{\beta_i}$
 is the preimage over $ Z^{min}_{\beta_i} = Z^{min}_{\alpha_i}$ under the
Bialynicki-Birula
contraction determined by $\beta_i$. If $\alpha_i$ and $\beta_i$ differ by a
{\sl positive}
constant multiple, we would have that $S^{l^+}_{\alpha_i}$ and
$S^{l^-}_{\beta_i}$ coincide.
Let us show that this is impossible. So,
assume that $S^{l^+}_{\alpha_i}$ and $S^{l^-}_{\beta_i}$ coincide.
Then pick  a point $z \in Z^{min}_{\alpha_i}$.
Since the map $f^-:X^s(l^-)/\!/G\to X^s(l_0)/\!/G$ is surjective, there is a
point $x \in X^s(l^-)$ such that
$G\cdot z \subset \overline{G\cdot x}$. Thus $x \in X^{sss}(l_0)$.
Since $x \in X^s(l^-) \cap X^{sss}(l_0)$, we have that $x \notin X^s(l^+)$
because $X^s(l^+) \cap X^s(l^-) = X^s(l_0)$ (Proposition 3.4.7).
Taking into account the stratification $X = X^s(l^+) \cup S^{l^+}_{\beta_1}
\cup \cdots \cup S^{l^+}_{\beta_q} \cup X^{us}(l_0),$
one sees that there must be $j \ne i$ such that $x \in S^{l^-}_{\beta_j}$. This
is because
we assume that $S^{l^+}_{\alpha_i}=S^{l^-}_{\beta_i}$. This shows that
there is a point $z' \in Z^{min}_{\beta_j} = Z^{min}_{\alpha_j}$ such that
$G\cdot z' \subset \overline{G\cdot x}$. This contradicts the fact that
$G\cdot z$ and $G\cdot z'$ are two distinct closed orbits in $X^{ss}(l_0)$.
Therefore $\alpha_i$ and $\beta_i$ differ by a {\sl negative}
constant multiple. This proves (ii) and (iii).

It remains to show (iv).
Assume that there was a point $x \in S^{l^-}_{\beta_j} \cap
S^{l^+}_{\alpha_i}$. Then
there must be a point $z \in Z^{min}_{\beta_j}$ and
a point $z' \in Z^{min}_{\alpha_i}$ such that
$\overline{G\cdot x} \supset G\cdot z$ and $\overline{G\cdot x} \supset G\cdot
z'$.
Because $x \in X^{ss}(l_0)$, the above would
imply that $G\cdot z$ and $G\cdot z'$ are mapped to the same point in the
quotient
$X^{ss}(l_0)/\!/G$. But this can not happen because  $G\cdot z$ and $G\cdot z'$
are different closed orbits in $X^{ss}(l_0)$.
\endproof

\smallskip
Let $\lambda_i$ be the one-parameter subgroup generated by $\beta_i$
and $\lambda_i^{-1}$ be the one-parameter subgroup generated by $\alpha_i$.
Since $F$ is a truly faithful cell, we obtain that
$G_z = \lambda_i({\bf C}^*)$ if $z \in Z^{min}_{\beta_i}$.
For convenience, we use the following notational convention.
For any $\beta\in \{\beta_1,\ldots,\beta_p\}$, we shall use $\alpha$ to denote
the corresponding element in $\{\alpha_1,\ldots,\alpha_p\}$ without specifying
the sub-index (cf. Lemma 4.2.4). The convention also extends to $\lambda$ and
$\lambda^{-1}$, and so on.

%\proclaim 4.2.5 Lemma. Keep the assumption and the notation from the previous
%lemma. Let $\beta\in \{\beta_1,\ldots,\beta_p\}.$ Then
%$$Y^{min}_{\beta} \setminus U(\lambda)  Z^{min}_{\beta} \subset X^s(l^-)$$
%and $$Y^{min}_{\alpha} \setminus U(\lambda^{-1})  Z^{min}_{\alpha} \subset
%%X^s(l^+).$$

%\proof Observe first that $Z^{min}_\beta = Z^{min}_{\alpha}$ and
%$Y^{min}_\beta \cap Y^{min}_{\alpha} = Z^{min}_\beta$ by the definition.

%Now since  $Y^{min}_\beta \subset X = X^s (l^-) \cup \bigcup_i S_{\alpha_i}
%%\cup X^{us}(l_0)$,
%$S_\beta \cap  X^{us}(l_0) = \emptyset$,  and
%$S_\beta \cap S_{\alpha_i} = \emptyset$ if $\beta \ne \beta_i$,
% we have that
%$Y^{min}_\beta \subset X^s (l^-) \cup S_{\alpha}$. Let $y \in
%%Y^{min}_{\beta}$.
%Then either $y \in X^s (l^-)$ or $y \in  S_{\alpha}$.
%Assume that $y \in S_{\alpha}$.
%Then there is a $g \in G$ and $y_1 \in Y^{min}_{\alpha}$
%such that $y = g\cdot y_1$. Because $G = U(\lambda) P(\lambda^{-1})$,
%we can write $g = u\cdot p $ with $u \in U(\lambda)$ and $ p \in
%%P(\lambda^{-1})$.
%Thus $y = up\cdot y_1 = u\cdot y_2$ where $y_2 = p\cdot y_1 \in
%%Y^{min}_{\beta}$.
%Thus $y_2 = u^{-1}\cdot y \in Y^{min}_\beta \cap Y^{min}_{\alpha} =
%%Z^{min}_\beta = Z^{min}_{\alpha}$
% (notice that $u \in U(\lambda)$ and
%$y \in Y^{min}_\beta$ implies $u^{-1}\cdot y \in Y^{min}_\beta$).
% Therefore $y = u\cdot y_2 \in U_\beta Z^{min}_\beta$. This implies that
%$$Y^{min}_{\beta} - U(\lambda)  Z^{min}_{\beta} \subset X^s(l^-).$$
%The other claim can be proved similarly.
%\endproof

\proclaim 4.2.5  Proposition. Keep the notation of 4.2.2 and the previous
assumption.
For any $\beta \in \{\beta_1, \cdots, \beta_p\}$, let
$$p_+: Y^{min}_{\beta} \rightarrow Z^{min}_{\beta}$$
$$p_-: Y^{min}_{\alpha} \rightarrow Z^{min}_{\alpha}$$
be the two natural projections.
The subgroup of $G$ that preserves the fiber of $p_\pm$ over $z \in
Z^{min}_{\beta}$
is $G_z \cdot U(\lambda)$ (resp. $G_z \cdot U(\lambda^{-1})$).

\proof  We shall consider only the map $p_+$. The other map is considered
similarly.
By 2.4.6 (i)  any element of $G$ which preserves the fibre
must belong to $P(\lambda)$. Let $p \in P(\lambda)$ and
$x \in p_+^{-1} (z)$ ($z \in Z^{min}_\beta$) such that $p \cdot x \in p_+^{-1}
(z)$.
Then we have $\lim_{t \to 0} \lambda(t) \cdot  x = z$ and $\lim_{t \to 0}
\lambda(t)
 p\cdot  x = z$.
But $\lim_{t \to 0} \lambda(t) \cdot p x =  \lim_{t \to 0} \lambda (t)
\cdot p \cdot \lambda^{-1}(t)
\lambda (t) x = p^\prime\cdot z$, where $p^\prime =\lim_{t \to 0} \lambda(t)
\cdot p \cdot \lambda^{-1}(t)$.

Hence $p^\prime z = z$ which  shows that $p^\prime \in G_z$. Hence
$\lim_{t \to 0} \lambda(t) \cdot  (p^\prime)^{-1} p \cdot \lambda^{-1}(t) =
{\rm id}$. That is,
$(p^\prime)^{-1} p \in U(\lambda)$, namely,
 $p \in p^\prime U(\lambda) \subset G_z U(\lambda)$, which implies the claim.
\endproof

\medskip
Following the above proposition, denote $p_\pm^{-1}(z)$ by $V^\pm$.
In what follows, we will concentrate on $V^+$. The other can be treated
similarly.

Now, the group $G_z$ acts on $U(\lambda)$ by conjugation.
Thus for any $u\cdot z\in U(\lambda)\cdot z$
and $g\in G_z$ we have $g\cdot (u\cdot z) = gug^{-1} g\cdot z = gug^{-1}\cdot
z$.
This shows that the
orbit $U(\lambda)\cdot z$ is $G_z$-invariant. Let us take a suitable
identification of $V^+$ with the affine space ${\bf C}^n$
so that the point $z$ is identified with the origin,  the group $G_z$
acts on $V^+$ linearly and has the point $z$ in the closure of any orbit. So
the action of
$G_z$ on $V^+$ is a good ${\bf C}^*$-action. It is well-known that it is
equivalent to
a positive
grading on the ring of regular functions ${\cal O}(V^+)\cong {\bf
C}[T_1,\ldots,T_n]$.
We assume that each
coordinate function $T_i$ is homogeneous of some degree
$q_i > 0 $. Let $R = U(\lambda)\cdot z$. Being an orbit of a unipotent group
acting on an affine
variety, $R$ is closed. Since it is $G_z$-invariant,
it can be given by a system of equations $F_1=\ldots =F_k = 0$ where $F_i$ are
(weighted) homogeneous generators of the ideal of regular functions on $V^+$
vanishing on $R$. We can write
$F_i = L_i(T_1,\ldots,T_n)+G_i(T_1,\ldots,T_n)$, where $L_i$ is a linear
function in $T_1,\ldots,T_n$ and
$G_i$ is a sum of (ordinary) homogeneous polynomials of degree $ > 1$. Since
$R$ is
nonsingular at the origin,
we can replace $F_i$ by
their linear combinations to assume that there exists $1\le s_1<\ldots <s_r \le
n, r = {\rm codim} \; R$, such that
$$L_i = T_{s_i}+\sum_{j\ne s_1,\ldots,s_k}a_{ij}T_j, i = 1,\ldots,r,\ \ L_i =
0, i> r.$$
Let $W$ be the linear subspace of $V^+$ defined by the equations
$T_j = 0, j\ne s_1,\ldots,s_r$. It is obviously $G_z$-invariant.
Consider the action map
$a:U(\lambda)\times W\to V^+$. It is $U(\lambda)G_z$-equivariant.

\proclaim 4.2.6 Lemma. The map $a$ is an isomorphism.

\proof
First we claim that $W\cap R = \{0\}$. In fact, since the equations $F_i$ are
weighted
homogeneous, each $G_i$ is a polynomial in $T_j, j\not\in \{s_1,\ldots,s_r\}$.
This shows that $R\cap W$ is given by the equations
$T_i = 0, i=1,\ldots,n$. This proves the claim. This also implies that the
tangent space at
the origin of
$V^+$ is equal to the direct sum of the tangent spaces of $R$ and $W$. Now the
source $V: = U(\lambda)\times W$ and the target $V^+$ are affine spaces of the
same dimension. The restriction of the map to $U(\lambda)\times \{0\}$ and to
$1\times W$ are
isomorphisms onto their images $R$ and $W$, respectively.
Thus the differential of the map $a$ at the point $(1, 0)$ is bijective.  Now
$$a^{-1}(a(1, 0)) = a^{-1}(z) = \{(u,w)\in U(\lambda)\times W : w = u^{-1}\cdot
z\}$$
consists of only one point $(1, z)$. We know that
the map is $G_z$-equivariant and the action of $G_z\cong {\bf C}^*$  on $V$ and
on
$V^+$ is a good ${\bf C}^*$-action. This implies that the map $a$ is defined by
a homomorphism of
positively graded
rings
$a^*:{\cal O}(V^+) \to {\cal O}(V)$. Let ${\bf m}_{V^+}$ be the maximal ideal
of the unique closed orbit
$\{(0,1)\}$ of $V$, and let ${\bf m}_V$ be the maximal ideal of the unique
closed orbit $\{z\}$ of $V^+$.
The property that $a$ is \'etale over $z$ and $a^{-1}(z) = (1,0)$ imply that
$a^*({\bf m}_{V^+}){\cal O}(V) = {\bf m}_V$. By
({\bf [Bo]}, Chapter III,\S1, Proposition 1), the algebra ${\cal O}(V)$ is
generated as
${\bf C}$-algebra by ${\bf m}_V$. This implies that the homomorphism $a^*$
is surjective. Since the both rings are
integral domains of the same
dimension, this shows that the map $a^*$ is an isomorphism. Hence $a$ is an
isomorphism.
\endproof

\proclaim 4.2.7 Corollary.
The quotient space of  $V^+  \setminus  U(\lambda) \cdot z$ by $G_zU(\lambda)$
can be identified with the quotient of $W  \setminus  \{0\}$ by $G_z$ which
is a weighted projective space. The similar statement holds for $V^-$.

\proof
First note that from Lemma 4.2.6 the orbit space of  $U(\lambda)$ on $V^+$
can be identified with $W$. In particular,
 the orbit space of $U(\lambda)$ on $V^+   \setminus  U(\lambda) \cdot z$
can be identified with $W  \setminus  \{z\}$.
This implies that the quotient of $V^+  \setminus  U(\lambda) \cdot z$ by
$G_zU(\lambda)$
can be identified with the quotient
of $W  \setminus  \{0\}$ by $G_z$ which is a weighted projective space
because the  $G_z$-action on $V^+$ (and hence on $W$) is good.
\endproof

\proclaim  4.2.8  Theorem.
Let $G$ be a reductive algebraic group acting on a nonsingular projective
variety $X$.
Let  ($C^+$,  $C^-$) be a pair of chambers relevant to a truly faithful cell
$F$.
Then,  there are two birational morphisms
$$f^+: X^s(C^+)/\!/G \rightarrow X^{ss}(F)/\!/G $$ and
$$f^-:  X^s(C^-)/\!/G \rightarrow X^{ss}(F)/\!/G $$
so that by
setting $\Sigma_0$ to be  $(X^{ss}(F) \setminus  X^s(F))/\!/G$,  we have
\item{(i)}  $f^+$ and $f^-$ are isomorphisms over the complement to $\Sigma_0$;
\item{(ii)} over each connected component $\Sigma_0'$ of  $\Sigma_0 $, the
fibres of the maps $f^\pm$ are
weighted projective spaces of dimension $d_{\pm}$ with the same weights;
\item{(iii)} $d_+ + d_-  + 1 = {\rm codim} \; \Sigma_0'$.

\proof
Statement (i) follows immediately from Lemma 4.2.1.

To show (ii) and (iii),  we apply the assertions and the notations from 4.2.2
and 4.2.6.

To prove (ii), it suffices to consider
any  particular stratum $S_\beta$. Let $\lambda$ be the corresponding
one-parameter
subgroup and
$$p_-: Y^{min}_{\alpha} \rightarrow Z^{min}_{\alpha},$$
$$p_+: Y^{min}_{\beta} \rightarrow Z^{min}_{\beta}$$
be the two natural projections.
Note that for any point $z\in Z^{min}_\alpha (= Z^{min}_\beta), G_z =
\lambda({\bf C}^*)$.

Let us denote the fiber $p^{-1}_{\pm} (z)$ by $V^\pm$.
The group $U(\lambda^\pm)$ acts on $V^\pm$
(not linearly!), $G_z$ acts on $V^+$ (resp. on $V^-$) linearly with positive
(resp. negative) weights.
We need only consider $V^+$; the other set can be treated similarly.

\smallskip
Now consider any non-stable
closed orbit $G \cdot z \subset X^{ss}(F)$ where $z \in Z^{min}_\alpha$
for some $\alpha$. We want to describe the fiber $(f^-)^{-1}([G \cdot z])$
(resp. $(f^+)^{-1}([G \cdot z])$) where $[G \cdot z] \in X^{ss}(F)/\!/G$ is the
induced
point in the quotient.
First we observe, by using Proposition 3.4.7, that
$$X^s(F) = X^s(C^+) \cap X^s(C^-).$$
By 4.2.2, we have
$$X = X^s(C^+) \cup S^{l^+}_{\alpha_1} \cup \cdots \cup S^{l^+}_{\alpha_p} \cup
X^{us}(F).$$
(Here $l^+ \in C^+$.) Thus by intersecting with $X^s(C^-)$ we obtain
$$X^s(C^-)= X^s(F) \cup  (X^s(C^-) \cap S^{l^+}_{\alpha_1}) \cup \cdots
 \cup  (X^s(C^-) \cap S^{l^+}_{\alpha_p}).$$
Now applying Lemma 4.2.3, we see that modulo the quotient relation in
$X^{ss}(F)$,
orbits from different  stratum $S^{l^+}_{\alpha_i}$
are identified with different closed orbits. Combining this with the above
decomposition
of $X^s(C^-)$, we conclude that
the orbits of $X^s(C^-)$ whose
induced points in the quotient $X^s(C^-)/G$
can be mapped to $[G \cdot z]$ by $f^-$ are contained in
$$X^s(C^-) \cap S^{l^+}_\alpha.$$
Next, observe that $S^{l^+}_\alpha = G Y^{min}_\alpha$.
So from Lemma 4.2.5, one sees that  the fiber of $(f^-)^{-1}([G \cdot z])$
can be identified with the orbit space
 of $X^s(C^-) \cap V^+$ by the group $G_z \cdot U(\lambda)$.
The set $X^s(C^-) \cap V^+$ is not empty because $f^-$ is surjective,
and it is also open in $V^+$ because $X^s(C^-)$ is open in $X$.
Furthermore, observe that
$$X^s(C^-) \cap V^+ \subset V^+  \setminus U(\lambda) \cdot z$$
 because any point
from $U(\lambda) \cdot z$ has isotropy subgroup of positive dimension. Hence
$$(X^s(C^-) \cap V^+) /G_z \cdot U(\lambda)$$ is open in
the quotient space of $V^+  \setminus U(\lambda) \cdot z$ by $G_zU(\lambda)$
(the latter is a weighted projective space by Corollary 4.2.7). But as the
fiber of $f^-$,
$X^s(C^-) \cap V^+ /G_z \cdot U(\lambda) $ is proper, thus we have that
$X^s(C^-) \cap V^+ /G_z \cdot U(\lambda) $ is equal to the weighted projective
space
$(V^+  \setminus U(\lambda) \cdot z)/G_zU(\lambda)$.
This completes the proof that the fiber of $f^-$ is a weighted projective
space.

\smallskip
Identical arguments
can be applied to obtain a similar result for the map $f^+$.
This finishes the proof of (ii).

\smallskip

To prove (iii),
we first claim that, for any $z \in Z^{min}_\alpha$, the set
$\{g \in G | g z \in Z^{min}_\alpha \}$ equals the normalizer
$N_\lambda$ of $G_z$.
To see this, if $g\cdot z \in Z^{min}_\alpha$ for some $z \in Z^{min}_\alpha$
then the isotropy
group $\lambda$ of $z$ should equal the
isotropy group $g \lambda g^{-1}$ of $g\cdot z$ because all of the elements of
$Z^{min}_\alpha$ have
 the same stabilizer $G_z$.
This implies that $g \in N_\lambda$. In fact, $Z^{min}_\alpha$ is just the set
of semi-stable
 points of the action of
$N^\prime_\lambda$ on $Z_\alpha$ where  $N^\prime_\lambda = N_\lambda / G_z$.

By the above claim and Theorem 2.4.6 (i),  we have
$$N_\lambda = \{g \in G | g z \in Z^{min}_\alpha \}
\subset \{g \in G | g z \in Y^{min}_\alpha \} = P(\lambda).$$ That is,
$N_\lambda$ is a subgroup in $P(\lambda)$. So any element $n$ of $N_\lambda$
can be
written as $ul$ with $u \in U(\lambda)$ and $l \in L(\lambda)$.
Now $L(\lambda) \subset N_\lambda$
because $L(\lambda)$ is
the centralizer of $\lambda$, so $u \in N_\lambda$.
That is, $u \in U(\lambda) \cap N_\lambda$.
But $N_\lambda$, as the normalizer of $G_z = \lambda({\bf C}^*)$,  is a
reductive subgroup of $G$. Therefore
$U(\lambda) \cap N_\lambda$ is finite.
This implies that $L(\lambda)$ is of finite index in $N_\lambda$.

To complete (iii),  let $\Sigma_0' = G Z^{min}_\alpha /\!/G = Z^{min}_\alpha
/\!/
 N^\prime_\lambda$.
Then by  the Bialynicki-Birula decomposition theorem ({\bf [B-B2]}), we have
$${\rm dim} \; Z^{min}_{\alpha} + {\rm dim} \; X =
{\rm dim} \; Y^{min}_{\alpha} +{\rm dim} \; Y^{min}_{\alpha}.$$
Therefore,
$${\rm dim} \; X -{\rm dim} \; Z^{min}_{\alpha} = {\rm dim} \; ({\rm fibre} \;
of \; p_+) +
 {\rm dim} \; ({\rm fibre} \; of \; p_-).$$
Hence
$${\rm dim} \; X -  {\rm dim} \; G  - ({\rm dim} \; Z^{min}_{\alpha} - {\rm
dim} \; N^\prime_\lambda)  $$
$$ = {\rm dim} \; ({\rm fibre} \; of \; p_+) - {\rm dim} \; U(\lambda)
 +  {\rm dim} \; ({\rm fibre} \; of \; p_-) - {\rm dim} \; U(\lambda^{-1}) -
1.$$
because ${\rm dim} \; G = {\rm dim} \; L(\lambda) +  {\rm dim} \; U(\lambda) +
{\rm dim} \; U(\lambda^{-1})$.
Now it follows that ${\rm codim} \; \Sigma_0' = d_+ + d_- + 1$.
\endproof

\smallskip\noindent
4.2.9 {\sl Remark.}  In the proof of (iii),
one really should treat connected components of $Z^{min}_{\alpha}$
separately since they may have different dimensions.
However this does not affect the proof at all (cf. 4.1.2).

\medskip\noindent
4.2.10 {\sl Remark.}  Theorem 4.2.8 can be generalized to the case of a
faithful (but not
necessary truly faithful) wall.
In this case,
 we can  still define the map $a:U(\lambda)\times W \to V^+$ which is
equivariant with respect to the one-dimensional radical of $G_z$.
We obtain that $W\setminus\{z\}$ is
contained in the union of $X^s(l^-)$ and some stratum $S_{\alpha_i}$. Thus the
orbit space
$V^+\cap X^s(l^-)/G_zU(\lambda)$ can be identified with the orbit space
$(W\setminus \{z\})\cap X^s(l^-)/G_z$. This agrees with an example showed to us
by C. Walters, where
the fibres are
isomorphic to Grassmannians and $G_z$ is isomorphic to $GL(n)$.

\eject
\smallskip
\bigskip\noindent
{\bf  References}.
\bigskip

\item{{\bf [At]}} M. Atiyah, Convexity and commuting Hamiltonians,
{\it Bull. London Math. Soc.} {\bf 14} (1982), 1-15.

\item{{\bf [Au]}} M. Audin, The topology of torus actions on symplectic
manifolds,
{\it Birkh\"auser}, 1991.

\item{{\bf [B-B1]}} A. Bialynicki-Birula, On homogeneous affine spaces of
linear algebraic groups, {\it
 Amer. J. Math.} {\bf 85} (1963), 577-582.

\item{{\bf [B-B2]}} A. Bialynicki-Birula, Some theorems on actions of algebraic
 groups,
{\it Ann. of Math.} {\bf 98} (1973), 480-497.

\item{{\bf [B-BS]}} A. Bialynicki-Birula and A.Sommese,
 Quotients by C* and SL(2,C) actions,
{\it Trans. Amer. Math. Soc.}
{\bf  279} (1983), 773-800.

\item {{\bf [Bo]}} N. Bourbaki,  {\it Commutative Algebra,}
Berlin, New York,  Springer-Verlag, 1989.

\item{{\bf [Br]}} M. Brion, Sur l'image de l'application moment,
in {\it Seminaire d'algebre Paul Dubreil et
Marie-Paule Malliavin, Paris, 1986}, {\it Lecture Notes in Math.} {\bf 1296},
177-193.

\item{{\bf [BP]}} M. Brion and C. Procesi,
 Action d'un tore dans une vari\'et\'e projective,
 in {\it ``Operator Algebras, Unitary Representations, Enveloping Algebras, and
Invariant Theory''},
     {\it   Progress in Mathematics, Birkh\"auser,}
 {\bf 192} (1990), 509-539.

%\item{{\bf [Do]}} I. Dolgachev, Introduction to Geometric Invariant Theory,
%%Lecture Notes
%Series, No 25, Seoul National University, 1994.

\item{{\bf [GM]}} M. Goresky and  R. MacPherson,
 On the topology of algebraic torus actions, in {\it ``Algebraic Groups,
Utrecht 1986''}
{\it Lect. Notes in Math.}
 {\bf 1271} (1986), 73-90.

%\item{{\bf [GH]}} P. Griffiths and J. Harris,
%{\it Principles of algebraic geometry}, John Wiley \& Sons, 1978.

%\item{{\bf [Gr]}} A. Grothendieck, Sur quelques points d'alg\`ebre
%%homologique,
%{\it Tohoku Math. J.} (2) {\bf 9} (1957), 119-221.

%\item{{\bf [GS1]}} V. Guillemin and  S. Sternberg,
%Convexity of the momentum mapping,
%{\it Invent. Math.}
%{\bf  67} (1982), 491-513.

\item{{\bf [GS]}} V. Guillemin and  S. Sternberg,
 Birational equivalence in the symplectic category,
{\it Invent. Math.}
{\bf  97} (1989), 485-522.

%\item{{\bf [GS3]}} V. Guillemin and  S. Sternberg,
%Geometric quantization and multiplicities of group representations,
%{\it Invent. Math.} {\bf 67} (1982), 515-538.

%\item{{\bf [Ha]}} R. Hartshorne,{\it Algebraic Geometry}, Springer.1977

\item{{\bf [He]}} W. Hesselink,
Desingularization of varieties of null forms,
{\it Inven. Math.} {\bf 55} (1979), 141-163.

\item{{\bf [Hu1]}} Y. Hu,
The geometric and topology of quotient varieties of torus actions,
{\it Duke Math. Journal}
 {\bf  68} (1992), 151-183.

\item{{\bf [Hu2]}} Y. Hu,
{\rm (W, R)} matroids and thin Schubert-type cells attached to algebraic torus
actions,
{\it Proc. of Amer. Math. Soc.} {\bf 123} No. 9 (1995), 2607-2617.

\item{{\bf [Ke]}} G. Kempf,  Instability in invariant theory,
{\it Ann. of Math.} {\bf 108} (1978), 299-316.

\item{{\bf [KN]}}  G. Kempf and L. Ness,
The length of vectors in representation spaces, in {\it ``Algebraic geometry,
Copenhagen 1978''},
{\it Lecture Notes in Math.} {\bf 732} (1979), Springer-Verlag,  pp.233-243.

\item{{\bf [Ki]}} F. Kirwan,
{\it Cohomology of quotients in symplectic and algebraic geometry},
Princeton University Press. 1984.

%\item{{\bf [Ki2]}} F. Kirwan,
%Convexity properties of the moment mapping, III
%{\it Invent. Math.} {\bf 77} (1984), 547-552.

\item{{\bf [KKV]}} F. Knop, H. Kraft, T. Vust,
 The Picard group of a G-variety,
in {\it ``Algebraic transformation groups and invariant theory''},
 DMV Seminar, B. 13, Birkha\"user , 77-87.

\item{{\bf [Kl]}} S. Kleiman, Towards a numerical criterion of ampleness,  {\it
Annals of Math.} (2) {\bf 84} (1966), 293-344.

\item{{\bf [KSZ]}} M. Kapranov, B. Sturmfels, A. Zelevinsky,
 Quotients of toric varieties,
{\it Math. Ann.} {\bf  290} (1991),
  643-655.

\item{{\bf [Li]}} D. Lieberman,
 Compactness of the Chow scheme:applications to automorphisms
and deformations of K\"ahler manifolds, in
{\it ``Seminair Francois Norguet 1975/77''},
{\it Lect. Notes in Math.} {\bf 670} (1978),  pp. 140-186.

\item{{\bf [MFK]}}  D. Mumford, J. Fogarty, F. Kirwan,
 Geometric Invariant Theory, 3d edition
{\it Springer-Verlag, Berlin, New York},
 1994.

\item{{\bf [Ne1]}} L. Ness, Mumford's numerical function and stable projective
hypersurfaces, in
{\it ``Algebraic geometry, Copenhagen 1978''}, {\it Lecture Notes in Math.}
{\bf 732} (1979), Springer-Verlag,
 pp. 417-453.

\item{{\bf [Ne2]}} L. Ness,
 A stratification of the null cone via the moment map,
{\it Amer. Jour.  of Math.},
 {\bf 106} (1984), 1281-1325.

%\item{{\bf [Ra]}} M. Raynaud, {\it Anneaux locaux hens\'eliens}, Lect. Notes
%%in Math. vol.169, Springer. 1970.

\item{{\bf [Re]}} M. Reid,
  What is a flip,
preprint, 1992,  17 pp.

\item{{\bf [Ro]}} R. Rockafellar, {\it Convex analysis}, Princeton Univ. Press,
1970.

\item{{\bf [Th1]}}  M. Thaddeus,
Stable pairs, linear systems and the Verlinde formula,
{\it Invent. Math.} {\bf 117} (1994), 317-353.

\item{{\bf [Th2]}}  M. Thaddeus,
Geometric invariant theory and flips, Journal of American Math. Society (to
appear).

\bigskip
\noindent Department of Mathematics, University of Michigan,
Ann Arbor, MI 48109

\medskip
\noindent {\it e-mail}:

\noindent
idolga@math.lsa.umich.edu

\noindent
yihu@math.lsa.umich.edu

\bigskip
\noindent {Current Address for} Y.H.:

\noindent
Department of Mathematics, University of Utah,
Salt Lake City, UT 84112

\medskip
\noindent {\it e-mail}:

\noindent
yhu@math.utah.edu

\end